\documentclass[a4paper,12pt,twoside]{article}
\usepackage{graphicx}\usepackage{a4wide}
\usepackage{amssymb,amsmath,amsthm}
\usepackage{siunitx}
\usepackage{subfigure,color,colordvi,graphicx,xcolor}
\usepackage[only,llbracket,rrbracket,interleave]{stmaryrd}
\usepackage[superscript,sort,compress]{cite}

\makeatletter
\newcommand*\wt[2][0.2ex]{%
        \begingroup
        \mathchoice{\wt@helper{#1}{#2}{\displaystyle}{\textfont}}
                   {\wt@helper{#1}{#2}{\textstyle}{\textfont}}
                   {\wt@helper{#1}{#2}{\scriptstyle}{\scriptfont}}
                   {\wt@helper{#1}{#2}{\scriptscriptstyle}{\scriptscriptfont}}%
        \endgroup
        #2%
}
\newcommand*\wt@helper[4]{%
        \def\currentfont{\the#41}%
        \def\currentskewchar{\char\the\skewchar\currentfont}%
        \setbox\tw@\hbox{\currentfont#2\currentskewchar}%
        \dimen@ii\wd\tw@
        \setbox\tw@\hbox{\currentfont#2{}\currentskewchar}%
        \advance\dimen@ii-\wd\tw@
        \rlap{\raisebox{-#1}{$\m@th#3\kern\dimen@ii\widetilde{\phantom{#2}}$}}%
}
\makeatother

\newenvironment{keywords}{\begin{quote}\emph{\textbf{Keywords:}}}{\end{quote}}

\newcommand{\bm}[1]{\text{\boldmath $#1$\unboldmath}}
\newcommand{\mat}[1]{\mathbf{#1}}
\newcommand{\nsd}{\ensuremath{\texttt{n}_{\texttt{sd}}}}
\newcommand{\msd}{\ensuremath{\texttt{m}_{\texttt{sd}}}}
\newcommand{\numel}{\ensuremath{\texttt{n}_{\texttt{el}}}}
\newcommand{\nface}{\ensuremath{\texttt{n}_{\texttt{fc}}}}
\newcommand{\grad}{\bm{\nabla}}
\newcommand{\gradS}{\bm{\nabla}_{\!\!\texttt{s}}}

\newcommand{\stressV}{\bm{\sigma}_{\!\!\texttt{v}}}

\newcommand{\bD}{\mat{D}}
\newcommand{\bDHalf}{\bD^{1/2}}

\newcommand{\bN}{\mat{N}}

\newcommand{\hu}{\hat{u}}
\newcommand{\hv}{\hat{v}}
\newcommand{\bn}{\bm{n}}

\newcommand{\bu}{\bm{u}}

\newcommand{\bw}{\bm{w}}
\newcommand{\bq}{\bm{q}}
\newcommand{\bhu}{\widehat{\bu}}
\newcommand{\bL}{\bm{L}}

\newcommand{\bt}{\bm{t}}
\newcommand{\Hdiv}{H(\operatorname{div})}
\newcommand{\Vh}{\ensuremath{\mathcal{V}^h \!}}
\newcommand{\HVh}{\ensuremath{\reallywidehat{\mathcal{V}}^h \!\!}}
\newcommand{\Poly}[1][1]{\ensuremath{\mathcal{P}^{#1}}}
\newcommand{\jump}[1]{\llbracket #1\rrbracket}
\newcommand{\Insd}{\mat{I}_{\nsd}}
\newcommand{\RR}{\mathbb{R}}
\newcommand{\sobo}[1][1]{\ensuremath{\mathcal{H}^{#1}}}
\newcommand{\eltwo}{\ensuremath{\mathcal{L}_2}}
\newcommand{\hDiv}[1]{\ensuremath{\mathcal{H}(\operatorname{div};{{#1}})}}
\usepackage{scalerel,stackengine}
\stackMath
\newcommand\reallywidehat[1]{%
\savestack{\tmpbox}{\stretchto{%
  \scaleto{%
    \scalerel*[\widthof{\ensuremath{#1}}]{\kern-.6pt\bigwedge\kern-.6pt}%
    {\rule[-\textheight/2]{1ex}{\textheight}}
  }{\textheight}%
}{0.5ex}}%
\stackon[1pt]{#1}{\tmpbox}%
}
\usepackage{xcolor}

\newtheorem{remark}{Remark}

\RequirePackage{algorithm}[0.1] 
\usepackage{algorithmic} 

\graphicspath{{figsPDF/}{figsPNG/}}

\begin{document}
\title{Discontinuous Galerkin approximations in computational mechanics: hybridization, exact geometry and degree adaptivity}

\author{
\renewcommand{\thefootnote}{\arabic{footnote}}
			  M. Giacomini\footnotemark[1]\textsuperscript{ \ ,}*, \
			  R. Sevilla\footnotemark[2]
}

\date{July 24, 2019}
\maketitle

\renewcommand{\thefootnote}{\arabic{footnote}}

\footnotetext[1]{Laboratori de C\`alcul Num\`eric (LaC\`aN), ETS de Ingenieros de Caminos, Canales y Puertos, Universitat Polit\`ecnica de Catalunya, Barcelona, Spain}
\footnotetext[2]{Zienkiewicz Centre for Computational Engineering, College of Engineering, Swansea University, Wales, UK
\vspace{5pt}\\
* Corresponding author: Matteo Giacomini. \textit{E-mail:} \texttt{matteo.giacomini@upc.edu}
}

\begin{abstract}
Discontinuous Galerkin (DG) discretizations with exact representation of the geometry and local polynomial degree adaptivity are revisited.
Hybridization techniques are employed to reduce the computational cost of DG approximations and devise the hybridizable discontinuous Galerkin (HDG) method.
Exact geometry described by non-uniform rational B-splines (NURBS) is integrated into HDG using the framework of the NURBS-enhanced finite element method (NEFEM).
Moreover, optimal convergence and superconvergence properties of HDG-Voigt formulation in presence of symmetric second-order tensors are exploited to construct inexpensive error indicators and drive degree adaptive procedures.
Applications involving the numerical simulation of problems in electrostatics, linear elasticity and incompressible viscous flows are presented.
Moreover, this is done for both high-order HDG approximations and the lowest-order framework of face-centered finite volumes (FCFV).
\end{abstract}


\section{Introduction}
\label{sc:Intro}

The importance of high-order approximations for the simulation of physical phenomena have been demonstrated in several fields of science and engineering, including electromagnetics~\cite{Hesthaven-HW-02,RS-DSM-18} and flow problems~\cite{Bassi-BR-97,AbgrallRicchiutoECM}.
DG methods have shown great potential for the development of efficient high-order discretizations, exploiting modern parallel computing architectures and adaptive strategies for non-uniform degree approximations~\cite{Karniadakis-CKS-00,Hesthaven-HW-02,Riviere2008,ErnBook,Cangiani2017}.
Nevertheless, the duplication of unknowns in classical DG methods and their resulting higher computational cost have limited their application mostly to academic problems and only few attempts to perform large-scale DG simulations are available in the literature, see~\cite{Bassi-CB-11,Munz-FM-16,Wall-FWK-19}.

To remedy this issue, \emph{static condensation} of finite element approximations~\cite{Guyan-65} and \emph{hybridization} of mixed methods~\cite{Fraeijs-65} have received special attention in recent years. 
Following the rationale in~\cite{Cockburn-16}, these concepts are applied to DG methods by defining the unknowns in each element as solution of a boundary value problem with Dirichlet data, whereas the interelement communication is handled by means of appropriate \emph{transmission conditions}.
Such approach leads to a wide range of \emph{hybrid discretization techniques}~\cite{brezzi1991mixed} in which the only globally-coupled degrees of freedom of the problem are located on the mesh faces.
The computational benefit of hybridization in the context of DG approximations has been analyzed in~\cite{AA-HARP:13} in terms of floating-point operations. Other thorough numerical comparisons are detailed in~\cite{Cockburn-KSC:11,May-WBMS-14}.

Contributions on hybrid methods may be subdivided in two main groups, relying either on primal or mixed formulations.
The former includes: (i) classical DG methods in which the number of coupled degrees of freedom is reduced simply by means of hybridization~\cite{Egger-ES-09,Egger-EW-12b,Egger-EW-12}; (ii) the \emph{reduced stabilization} approach exploiting a primal unknown approximated using a polynomial function of degree $k {+} 1$ and a trace variable of polynomial degree $k$ to furtherly ease the computational burden~\cite{Oikawa-15,Oikawa-16}; (iii) the \emph{hybrid high-order} (HHO) method which introduces a local reconstruction operator to mimick the behavior of the gradient of the primal solution and an appropriate stabilization term~\cite{DiPietro-DPEL-14,Ern-DPE-15}.
It is worth recalling that the HHO method belongs to the family of hybridizable DG approaches and can be recasted in this framework by an appropriate definition of the involved stabilization operator~\cite{Ern-CDPE-16}.

Stemming from the work on the local DG method~\cite{Cockburn-CS-98,Cockburn-CDG:08}, the hybridizable DG method proposed by Cockburn and coworkers relies on a mixed hybrid formulation~\cite{brezzi1991mixed}, based on polynomial approximations discontinuous element-by-element~\cite{cockburn2004characterization}.
The latter group thus includes all HDG formulations featuring the introduction of a mixed variable~\cite{Jay-CG-09,Jay-CGL:09,Nguyen-NPC:09,Nguyen-NPC:09b,Nguyen-CNP:10,Nguyen-NPC:10,Nguyen-NPC:11,soon2009hybridizable,Fu-FCS-15,May-SM-13}.
The advantage of directly approximating flux/stress via the introduction of a mixed variable is of special interest in the context of engineering problems in which quantities of interest usually rely on such information.
Thus, in the following sections, these specific hybrid methods based on mixed formulations will be considered and, with an abuse of notation, they will be denoted generically as HDG approaches.

Hybrid discretization techniques have been successfully applied to several problems of engineering interest.
In the context of computational fluid dynamics, HDG mixed formulations of the incompressible Navier-Stokes equations have been presented in~\cite{Nguyen-NPC:11,Cesmelioglu-CCQ-17} and~\cite{Shi-QS-16} using equal order and different order of polynomial approximations for the primal, mixed and hybrid variables, respectively. HHO formulations have been discussed in~\cite{DiPietro-DPK-18,DiPietro-BDPD-19}.
On the one hand, special emphasis has been devoted to the construction of pointwise divergence-free approximations in incompressible flows~\cite{Lehrenfeld-LS-16,Wells-RW-18b}. Recent results proposing a relaxed $\Hdiv$-conforming discretization of the velocity field are available in~\cite{Lehrenfeld-LLS-18,Lehrenfeld-LLS-19}.
On the other hand, extension to turbulent flows using implicit large eddy simulations~\cite{Peraire-FNP-17} and the Spalart-Allmaras model~\cite{Peraire-MNP-11,Evans-PE-19} and treatment of complex rheologies like quasi-Newtonian fluids~\cite{Gatica-GS-15} and viscoplastic materials~\cite{Ern-CBCE-18} are active topics of investigations.
First results of the application of hybrid discretization techniques to compressible flow problems are available in~\cite{peraire2010hybridizable,Williams-18}.

Concerning linear elasticity, the strong enforcement of the symmetry of the stress tensor in HDG has been studied by different authors. A formulation using different degrees of polynomial approximation for the primal and hybrid variables has been discussed in~\cite{Shi-QSS-18}. In~\cite{Cockburn-CF-17}, an appropriate enrichment of the local discrete space of approximation via the $M$-decomposition framework is proposed to ensure optimal convergence of the mixed variable and superconvergence of the postprocessed one. An easy-to-implement alternative is represented by the HDG-Voigt approach introduced in~\cite{RS-SGKH:18} and detailed in Section~\ref{sc:HDG-Voigt} of the present contribution.
In the context of nonlinear elasticity, hybrid methods based on primal formulations have shown promising results, see~\cite{Ern-AEP-18,Ern-AEP-19,Ern-AEP-19b} for HHO applications to hyperelastic, plastic and elastoplastic regimes. The exploitation of HDG mixed formulations to simulate these phenomena is currently an open problem, as described in~\cite{Cockburn-KLC-15,Cockburn-CS-19,Peraire-TNBP-19}.
Moreover, results on fluid-structure interaction problems and arbitrary Lagrangian Eulerian formulations have been investigated in~\cite{Pitt-SMP-16} and~\cite{Fidkowski-16}, respectively.

Other fields actively studied using hybrid discretization methods include subsurface flows~\cite{Riviere-FKR-18,Sarrate-CRS-19} and wave propagation phenomena~\cite{Peraire-FCTNP-18}, spanning from elastodynamics~\cite{Lanteri-BCDL-17,Terrana-TVG-17,Sayas-HPS-17} to coastal water simulations~\cite{Dawson-SD-18}, from Maxwell's equations~\cite{Lanteri-CDL-18} to acoustics~\cite{Wall-SKW-18}, optics~\cite{Lanteri-LLMW-17} and plasmonics~\cite{Peraire-VCNOP-18,Peraire-VCNP-18}.

Besides the application of hybrid discretization methods to different physical problems, several efforts have been devoted in recent years to the construction of efficient strategies to exploit the numerical advantages of the above mentioned approaches.
On the one hand, the flexibility of DG methods has been exploited to perform mesh refinement based on octrees~\cite{Dawson-SMD-16}, driven by adjoint-based~\cite{May-WMS-14} and fully-computable~\cite{Fu-AF-18} \emph{a posteriori} error estimators. On the other hand, the possibility of using nonuniform polynomial degree approximations has been explored in~\cite{Wall-HBKPCW-18} and~\cite{RS-SH:18,RS-19}. 

It is worth recalling that the accuracy of the functional approximation is strictly related to the one of the geometrical description of the domain. In this context, HDG for domains with curved boundaries have been analyzed in~\cite{Solano-CS-12,Solano-CS-14,Solano-SV-19,SanchezVizuet-SVS-19} via the extension to a fictitious subdomain, whereas a classical isoparametric framework has been developed for HHO in~\cite{DiPietro-BDP-18}. In~\cite{RS-SH:18,RS-19}, the NEFEM paradigm is  coupled with HDG to treat exact geometries described by means of NURBS. The strict relationship between geometrical and functional approximation error and its importance in the context of degree-adaptive procedures is further detailed in Section~\ref{sc:HDG-NEFEM} of the present contribution.

Recently, different approaches to problems featuring unfitted interfaces have been proposed using immersed HDG formulations~\cite{unfittedHDG,Solano-QSV-16}, the extended HDG framework (X-HDG) which mutuates ideas from X-FEM to treat cut cells~\cite{Gurkan-GSKF-16,Gurkan-GKF-17,Gurkan-GKF-19} and the cut-HHO method which relies on a cell agglomeration procedure and exploits the capability of HHO to handle generic mesh elements~\cite{Ern-EB-18}.
Moreover, numerical strategies to couple continuous Galerkin and HDG discretizations have also been recently proposed for mono- and multiphysics problems~\cite{Paipuri-PTF-19,LaSpina-LSGH-19}.

Concerning specific solution strategies for hybrid discretization methods, a parallel solver based on the iterative Schwarz method has been developed in~\cite{Gander-GH-18}, fast multigrid solvers have been employed in~\cite{Schutz-SA-17,Wall-KW-18,Riviere-FKMR-19} and iterative approaches inspired by the Gauss-Seidel method have been discussed in~\cite{BuiThanh-MTBT-17,BuiThanh-MTBT-18}.
Moreover, tailored preconditioners for the hybrid DG method have been proposed in~\cite{Wells-RW-18,Dolean-BBDNT-18} in the context of the Stokes equations.

This contribution presents an overview of some recent advances on HDG methods with application to different problems in computational mechanics, namely electrostatics, linear elasticity and incompressible viscous flow simulations.
The rationale to devise an HDG mixed approximation of a second-order partial differential equation (PDE) is recalled in Section~\ref{sc:HDG}.
In Section~\ref{sc:HDG-NEFEM}, the importance of accounting for the exact geometry described by means of NURBS is illustrated via the framework of NEFEM.
An HDG-NEFEM discretization with degree adaptivity is thus discussed for an electrostatics problem.
In Section~\ref{sc:HDG-Voigt}, an application of HDG to linear elasticity is considered.
Special attention is devoted to the construction of a formulation using a pointwise symmetric mixed variable, namely the strain rate tensor, via Voigt notation~\cite{FishBelytschko2007}.
The resulting HDG-Voigt formulation is robust for nearly-incompressible materials and provides optimally-convergent stresses and superconvergent displacements which are exploited to construct local error indicators to perform degree adaptive procedures.
Eventually, a lowest-order HDG approximation, the recently proposed FCFV method \cite{RS-SGH:18,RS-SGH:19}, is devised to efficiently solve large-scale problems involving incompressible flows (Section \ref{sc:FCFV}).
The FCFV method provides an LBB-stable discretization which is insensitive to mesh distortion and stretching and features first-order accurate fluxes without the need to perform a reconstruction procedure.

\section{The HDG rationale}
\label{sc:HDG}

To recall the rationale of the HDG method, the Laplace equation is considered in an open bounded domain $\Omega \subset \RR^{\nsd}$, $\nsd$ being the number of spatial dimensions, 
\begin{equation} \label{eq:Poisson}
\left\{\begin{aligned}
-\grad{\cdot}\grad u &= 0       &&\text{in $\Omega$,}\\
u &= u_D  &&\text{on $\partial \Omega$,}\\
\end{aligned}\right.
\end{equation}
where $u$ and $u_D$ are the unknown variable and its imposed value on the boundary, respectively.
From the point of view of modeling, Equation~\eqref{eq:Poisson} represents an electrostatic problem where $u$ is the unknown electric potential.

The standard HDG mixed formulation described in~\cite{RS-SH:16} is detailed. 
Recall that the main features of this HDG method is the introduction of a mixed variable, namely $\bq {=} {-}\grad u$ allowing to rewrite a second-order PDE as a system of first-order PDEs, and of a hybrid variable $\hu$ representing the trace off the primal unknown on the faces of the \emph{internal skeleton}
\begin{equation*}
\Gamma := \left[ \bigcup_{e=1}^{\numel} \partial\Omega_e \right] \setminus \partial\Omega \, ,
\end{equation*}
where $\numel$ is the number of non-overlapping elements $\Omega_e, \, e{=}1,\ldots,\numel$ in which the domain is partitioned.
Thus, Equation~\eqref{eq:Poisson} is rewritten as a system of first-order PDEs element-by-element 
\begin{equation*}
\left\{\begin{aligned}
\bq_e + \grad u_e &= \bm{0} &&\text{in $\Omega_e, \, e=1,\ldots,\numel$,}\\
\grad {\cdot} \bq_e &= 0  &&\text{in $\Omega_e, \, e=1,\ldots,\numel$,}\\
u_e &= u_D  &&\text{on $\partial\Omega_e\cap\partial\Omega$,}\\
u_e &= \hu  &&\text{on $\partial\Omega_e\setminus\partial\Omega$,}\\
\end{aligned}\right.
\end{equation*}
with the following \emph{transmission conditions} enforcing the continuity of the solution and of the fluxes across the interface $\Gamma$
\begin{equation*}
\left\{\begin{aligned}
\jump{u\bn} &= \bm{0} &&\text{on $\Gamma$,}\\
\jump{\reallywidehat{\bn {\cdot} \bq}} &= 0  &&\text{on $\Gamma$,}
\end{aligned}\right.
\end{equation*}
where $\jump{\odot} = \odot_i + \odot_l$ is the \emph{jump} operator proposed in~\cite{AdM-MFH:08} as the sum of the values in the elements $\Omega_i$ and $\Omega_l$ on the right and on the left of the interface respectively, whereas the trace of the numerical flux is defined as
\begin{equation*}
\reallywidehat{\bn {\cdot} \bq} := 
\begin{cases}
\bn {\cdot} \bq_e + \tau (u_e - u_D) & \text{on $\partial\Omega_e\cap\partial\Omega$,} \\
\bn {\cdot} \bq_e + \tau (u_e - \hu) & \text{elsewhere,}  
\end{cases}
\end{equation*}
with $\tau$ being an appropriate stabilization parameter~\cite{Jay-CGL:09,Nguyen-NPC:09,Nguyen-NPC:09b,Nguyen-CNP:10,Nguyen-NPC:11}.
Note that the first transmission condition is automatically fulfilled owing to the Dirichlet boundary condition $u_e{=}\hu$ imposed in the local problems on $\partial\Omega_e\setminus\partial\Omega$ and to the uniqueness of the hybrid variable $\hu$ on each mesh face in $\partial\Omega_e \subset \Gamma$.

Thus, the HDG local problems are defined as follows: for $e{=}1,\ldots,\numel$ compute $(u_e,\bq_e) \in \sobo(\Omega_e) \times \left[\hDiv{\Omega_e};\RR^{\nsd}\right]$ such that
%
\begin{equation} \label{eq:HDGPoissonLoc}
\left\{\begin{aligned}
- \int_{\Omega_e}{\!\!\! \bw {\cdot} \bq_e \, d\Omega} + \int_{\Omega_e}{\!\!\! \grad {\cdot} \bw \, u_e \, d\Omega} &= \int_{\partial\Omega_e\cap\partial\Omega}{\!\!\! \bn {\cdot} \bw \, u_D \, d\Gamma} + \int_{\partial\Omega_e\setminus\partial\Omega}{\!\!\! \bn {\cdot} \bw \, \hu \, d\Gamma} ,\\
\int_{\Omega_e}{v \, \grad {\cdot} \bq_e \, d\Omega} + \int_{\partial\Omega_e}{\!\!\! \tau \, v \, u_e \, d\Gamma} &= \int_{\partial\Omega_e\cap\partial\Omega}{\!\!\! \tau \, v \, u_D \, d\Gamma} + \int_{\partial\Omega_e\setminus\partial\Omega}{\!\!\! \tau \, v \, \hu \, d\Gamma} ,
\end{aligned}\right.
\end{equation}
for all $(v,\bw) \in \sobo(\Omega_e) \times \left[\hDiv{\Omega_e};\RR^{\nsd}\right]$, where $\left[\hDiv{\Omega_e};\RR^{\nsd}\right]$ is the space of square integrable vectors of dimension $\nsd$ with square integrable divergence on $\Omega_e$.

Following the notation in~\cite{RS-SH:16}, the discrete functional spaces
\begin{subequations}\label{eq:HDGspaces}
\begin{align}
\Vh(\Omega) & {:=}
\left\{ 
v \in \eltwo(\Omega) : \, v\vert_{\Omega_e} \in \Poly[k](\Omega_e) \, \forall \Omega_e, \, e {=} 1,\ldots,\numel
\right\},
\\
\HVh(S) & {:=}
\left\{
\hv \in \eltwo(S) : \, \hv\vert_{\Gamma_i} \in \Poly[k](\Gamma_i) \, \forall \Gamma_i \subset S \subseteq \Gamma \cup \partial\Omega
\right\} ,
\end{align}
\end{subequations}
are introduced for the HDG approximation. In Equation~\eqref{eq:HDGspaces}, $\Poly[k](\Omega_e)$ (respectively, $\Poly[k](\Gamma_i)$) represents the space of polynomial functions of complete degree at most $k {\geq} 1$ in $\Omega_e$ (respectively, on $\Gamma_i$).
Thus, for $e {=} 1,\ldots,\numel$ the HDG discrete local problem is: given $\hu^h$ on $\Gamma$, find $(u_e^h,\bq_e^h) \in \Vh(\Omega_e) {\times} \left[\Vh(\Omega_e)\right]^{\nsd}$, approximating the pair $(u_e,\bq_e)$, such that Equation~\eqref{eq:HDGPoissonLoc} holds for all $(v,\bw) \in \Vh(\Omega_e) {\times} \left[\Vh(\Omega_e)\right]^{\nsd}$.

\begin{remark}
For each element $\Omega_e, \, e{=}1,\dots,\numel$, the primal, $u_e^h$, and mixed, $\bq_e^h$, variables are determined as functions of the unknown hybrid variable $\hu^h$ on $\partial\Omega_e\setminus\partial\Omega$. 
From the point of view of modeling, the HDG local problem establishes a relationship between the electric potential and electric field inside each element and the electric potential on the corresponding element boundary.
\end{remark}

The HDG global problem is defined from the previously introduced transmission conditions: find $\hu \in \sobo[1/2](\Gamma)$ such that
%
\begin{equation} \label{eq:HDGPoissonGlob}
\sum_{e=1}^{\numel} \left\{ 
\int_{\partial\Omega_e\setminus\partial\Omega}{\!\!\! \hv \, \bn {\cdot} \bq_e \, d\Gamma} + \int_{\partial\Omega_e\setminus\partial\Omega}{\!\!\! \tau \, \hv \, u_e \, d\Gamma} - \int_{\partial\Omega_e\setminus\partial\Omega}{\!\!\! \tau \, \hv \, \hu \, d\Gamma} \right\} = 0 ,
\end{equation}
for all $\hv \in \eltwo(\Gamma)$, where $u_e$ and $\bq_e$ are obtained from the local problems defined in Equation~\eqref{eq:HDGPoissonLoc}.

The HDG discrete global problem is thus obtained solving the previous equation in the hybrid space introduced in~\eqref{eq:HDGspaces}, that is, find $\hu^h \in \HVh(\Gamma)$ such that Equation~\eqref{eq:HDGPoissonGlob} holds for all $\hv \in \HVh(\Gamma)$.

Recall that using equal order $k$ for the approximation of the primal, mixed and hybrid variables, HDG provides optimal convergence of order $k{+}1$ for all the unknowns~\cite{Jay-CGL:09}. Inspired by the work of Stenberg~\cite{Stenberg-88}, this property is exploited to devise an inexpensive local postprocessing procedure leading to a superconvergent approximation of the primal variable~\cite{Nguyen-NPC:10,RS-SH:16,RS-SH:18}: for $e{=}1,\ldots,\numel$, compute $u_e^\star$ using a polynomial approximation of degree $k{+}1$ such that
\begin{equation} \label{eq:HDGPoissonPost}
\left\{\begin{aligned}
\grad{\cdot}\grad u_e^\star &= -\grad{\cdot}\bq_e^h  &&\text{in $\Omega_e$,}\\
\bn {\cdot} \grad u_e^\star &= -\bn {\cdot} \bq_e^h  &&\text{on $\partial\Omega_e$,}\\
\end{aligned}\right.
\end{equation}
with the solvability constraint
\begin{equation}\label{eq:PostProcessCondMean}
\int_{\Omega_e}{\!\!\! u_e^\star \, d\Omega} = \int_{\Omega_e}{\!\!\! u_e^h \, d\Omega} .
\end{equation}
The computed $u_e^\star$ thus superconverges with order $k{+}2$~\cite{cockburn2008superconvergent} and has been exploited to define a simple and inexpensive error indicator to perform degree adaptive procedures~\cite{giorgiani2013hybridizable,giorgiani2014hybridizable,RS-SH:18}.

Henceforth, the superscript $^h$ identifying the discrete HDG solution will be omitted to ease readability and notation, if no risk of ambiguity is possible.

\section{HDG-NEFEM: exact geometry and degree adaptivity}
\label{sc:HDG-NEFEM}

The possibility to easily implement a variable degree of approximation in DG methods has motivated the recent interest in degree adaptive processes for convection-dominated flow and wave propagation phenomena. 
In this context, the superconvergent property of HDG is especially attractive, as it allows to devise an inexpensive error indicator for a computed approximation~\cite{giorgiani2013hybridizable,giorgiani2014hybridizable,RS-SH:18}.

One aspect that has been traditionally ignored when proposing new degree adaptive procedures is the representation of the geometry. 
In an isoparametric formulation, a degree adaptive process requires communicating with the CAD model and regenerating the mesh, at least near the boundary, at each iteration.
Nonetheless, the associated computational cost makes this strategy unfeasible for practical applications.
Thus, it is common practice to represent the geometry with quadratic or cubic polynomials and change only the degree of the functional approximation during the adaptivity process, leading to subparametric and superparametric formulations~\cite{giorgiani2013hybridizable,giorgiani2014hybridizable}. 
An alternative procedure based on the NEFEM rationale~\cite{RSC-SFH:08} is discussed here. 
The boundary of the computational domain is represented using the true CAD model, irrespective of the functional approximation used. 
The effort required to implement this approach is similar to the one employing subparametric or superparametric formulations: no communication with the CAD model or regeneration of the mesh are required, while the geometric uncertainty introduced by a polynomial description of the boundary of the domain is completely removed.

According to the framework described in Section~\ref{sc:HDG}, for a given distribution of the degree $k$ of the functional approximation, the global problem~\eqref{eq:HDGPoissonGlob} is solved first to obtain the trace of the electric potential on the mesh edges/faces. 
Second, an element-by-element problem is solved to compute the value of the electric potential $u$ and its gradient, i.e. the electric field, in the elements, according to Equation~\eqref{eq:HDGPoissonLoc}. 
Finally, an element-by-element postprocess is performed to obtain a superconvergent solution $u_e^\star$ by solving \eqref{eq:HDGPoissonPost}-\eqref{eq:PostProcessCondMean}. 

Following~\cite{giorgiani2014hybridizable,RS-SH:18}, a measure of the error in each element $\Omega_e, e=1,\ldots,\numel$ is defined
\begin{equation} \label{eq:errorMeasureU}
E_e^u = \left[ \frac{1}{|\Omega_e|} \int_{\Omega_e} \left( u_e^\star - u_e \right)^2 d\Omega \right]^{1/2}.
\end{equation}
Moreover, the local \emph{a priori} error estimate derived in~\cite{diez1999unified} for elliptic problems, states that the error in an element is bounded as
\begin{equation} \label{eq:aPrioriError}
\varepsilon_e = \| u - u_h \|_{\Omega_e} \leq C h_e^{k_e + 1 + \nsd/2} .
\end{equation}
By means of Richardson extrapolation, it is possible to estimate the unknown constant $C$ in Equation~\eqref{eq:aPrioriError}, assuming that two values of the error, obtained with different degrees of approximation, are considered. 
In order to determine the change of degree required to achieve a desired error $\varepsilon$, first, an estimate of the error is devised using Equation~\eqref{eq:errorMeasureU}. Then, the target approximation degree is computed according to Equation~\eqref{eq:aPrioriError} by imposing a desired elemental error $\varepsilon_e$. As detailed in~\cite{RS-SH:18}, the change of degree in the element $\Omega_e$ is thus given by
\begin{equation} \label{eq:newP}
\Delta k_e = \left\lceil \frac{\log(\varepsilon/E_e^u)}{\log(h_e)} \right\rceil, \quad e=1,\ldots,\numel ,
\end{equation}
where $\lceil \cdot \rceil$ is the ceiling function and $h_e$ is the non-dimensional characteristic size of the element $\Omega_e$.

\begin{figure}
	\centering
	\includegraphics[width=0.3\textwidth]{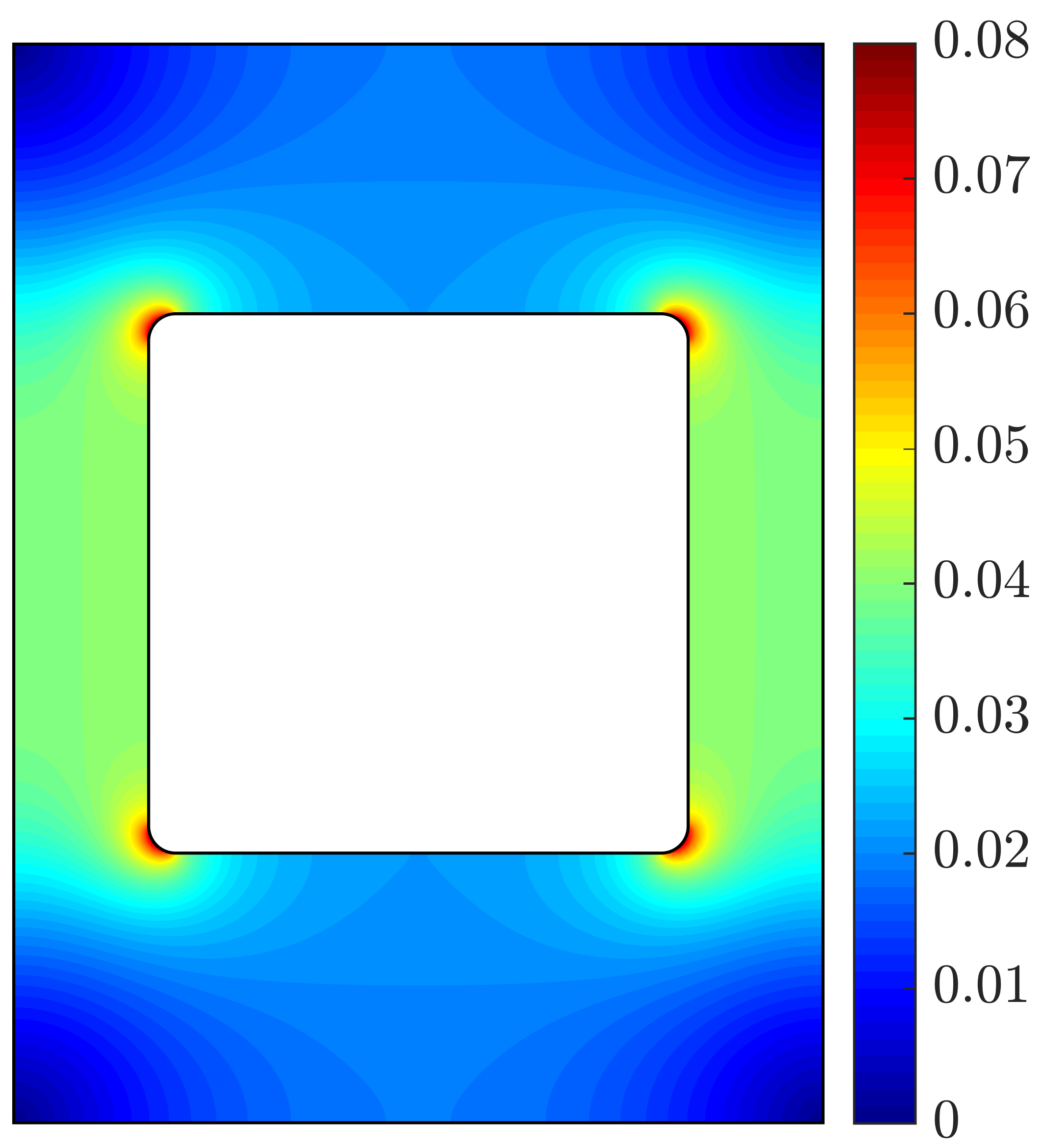}
	\includegraphics[width=0.3\textwidth]{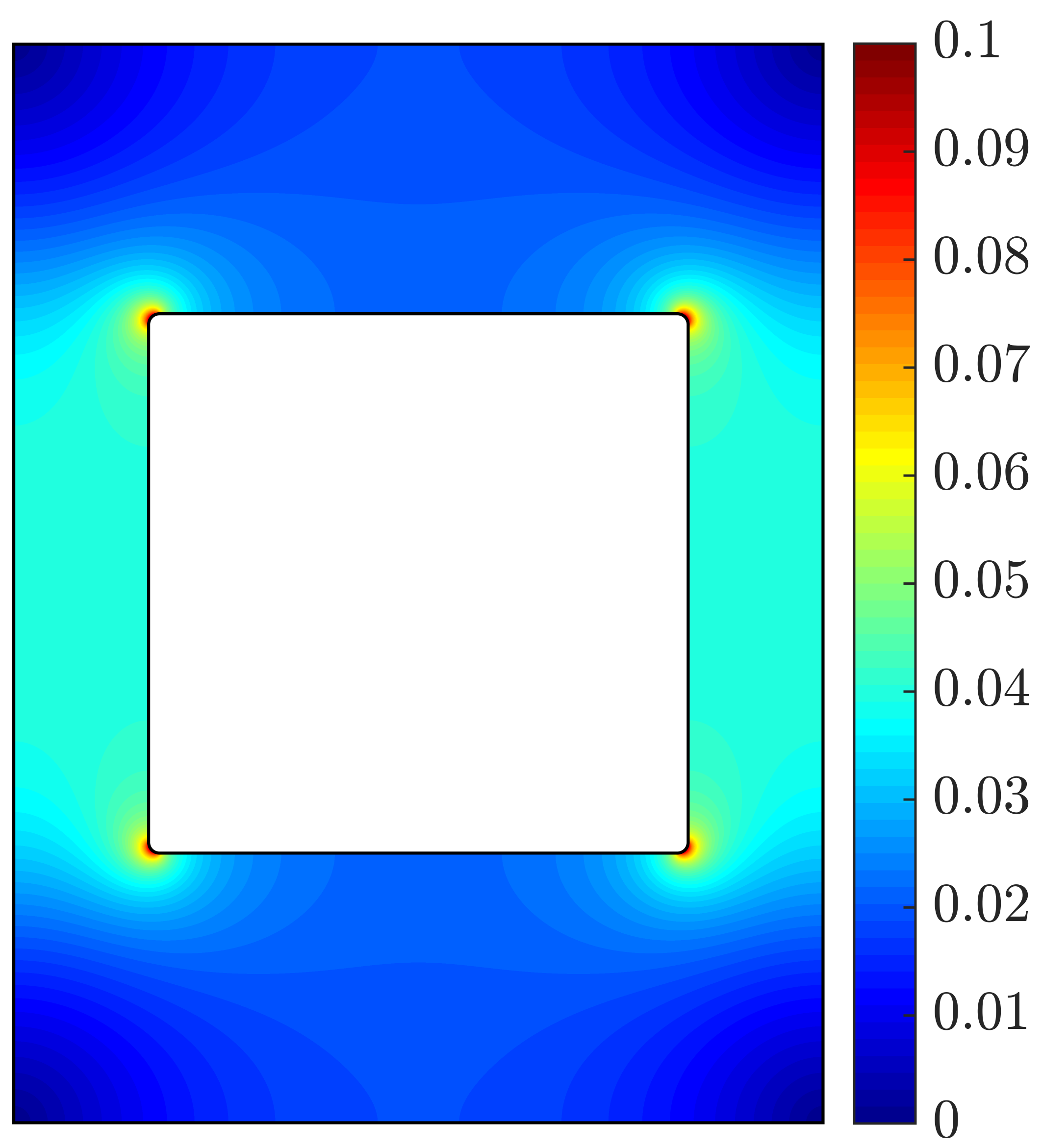}
	\includegraphics[width=0.3\textwidth]{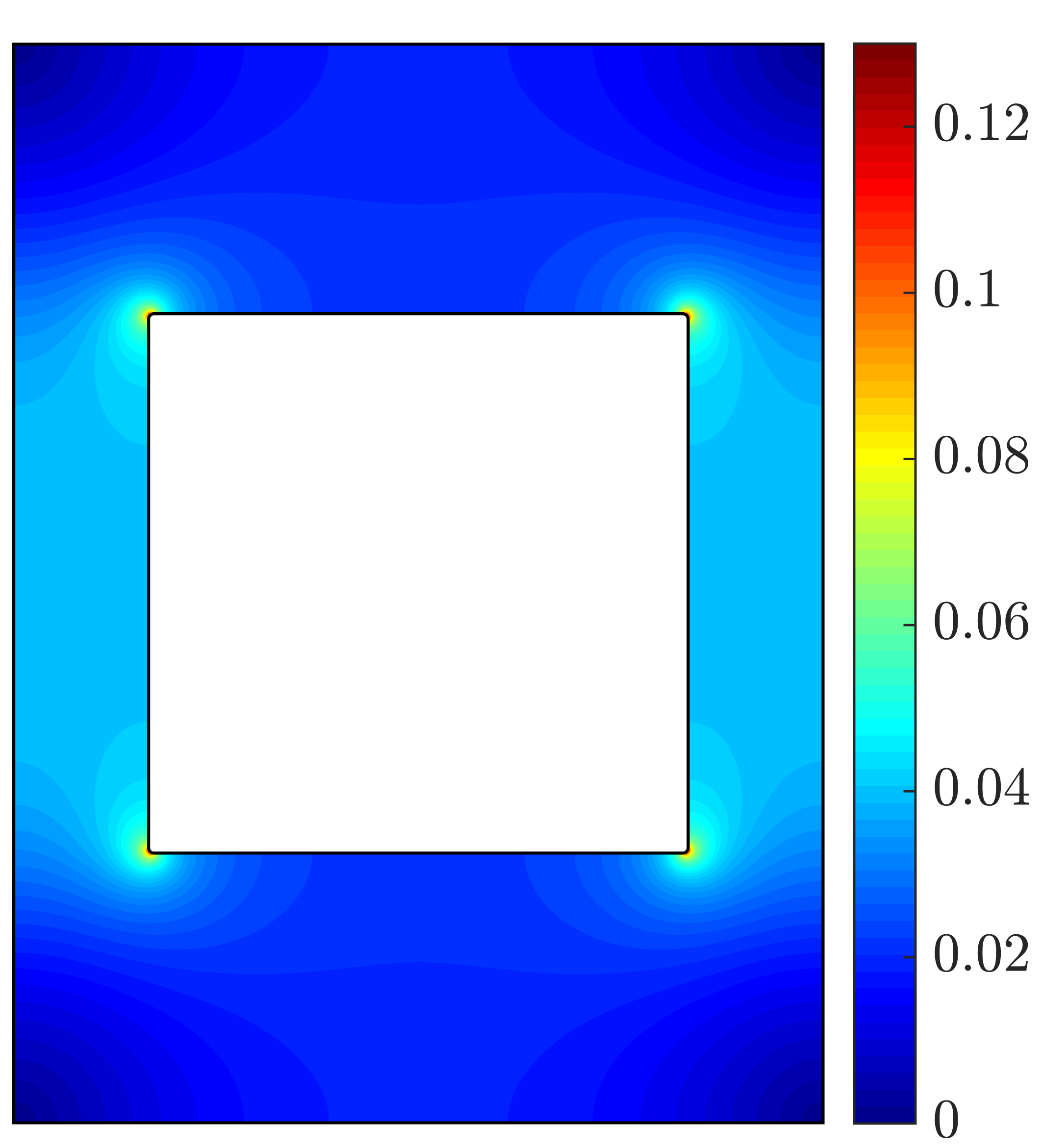}
	
	\includegraphics[width=0.3\textwidth]{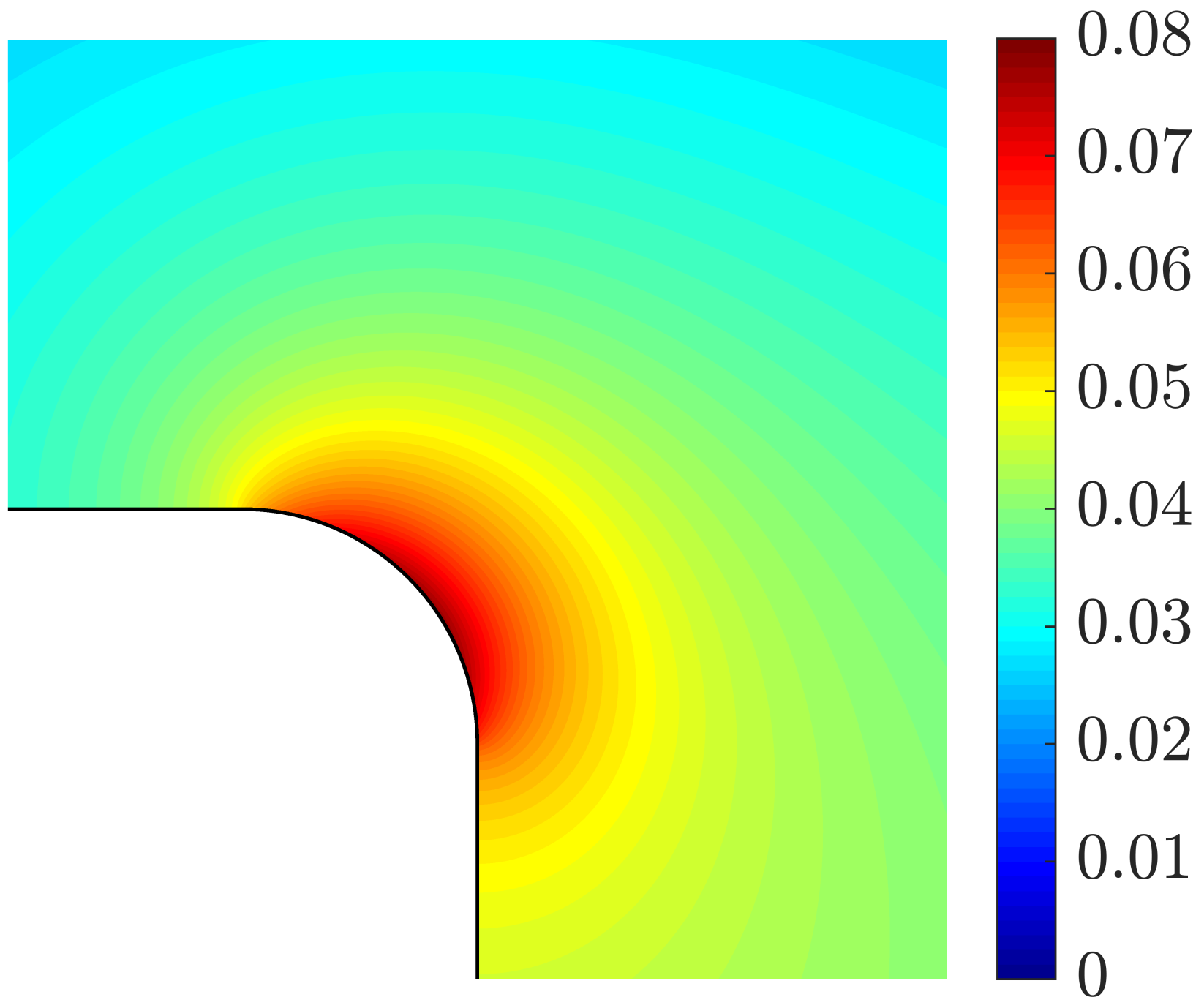}
	\includegraphics[width=0.3\textwidth]{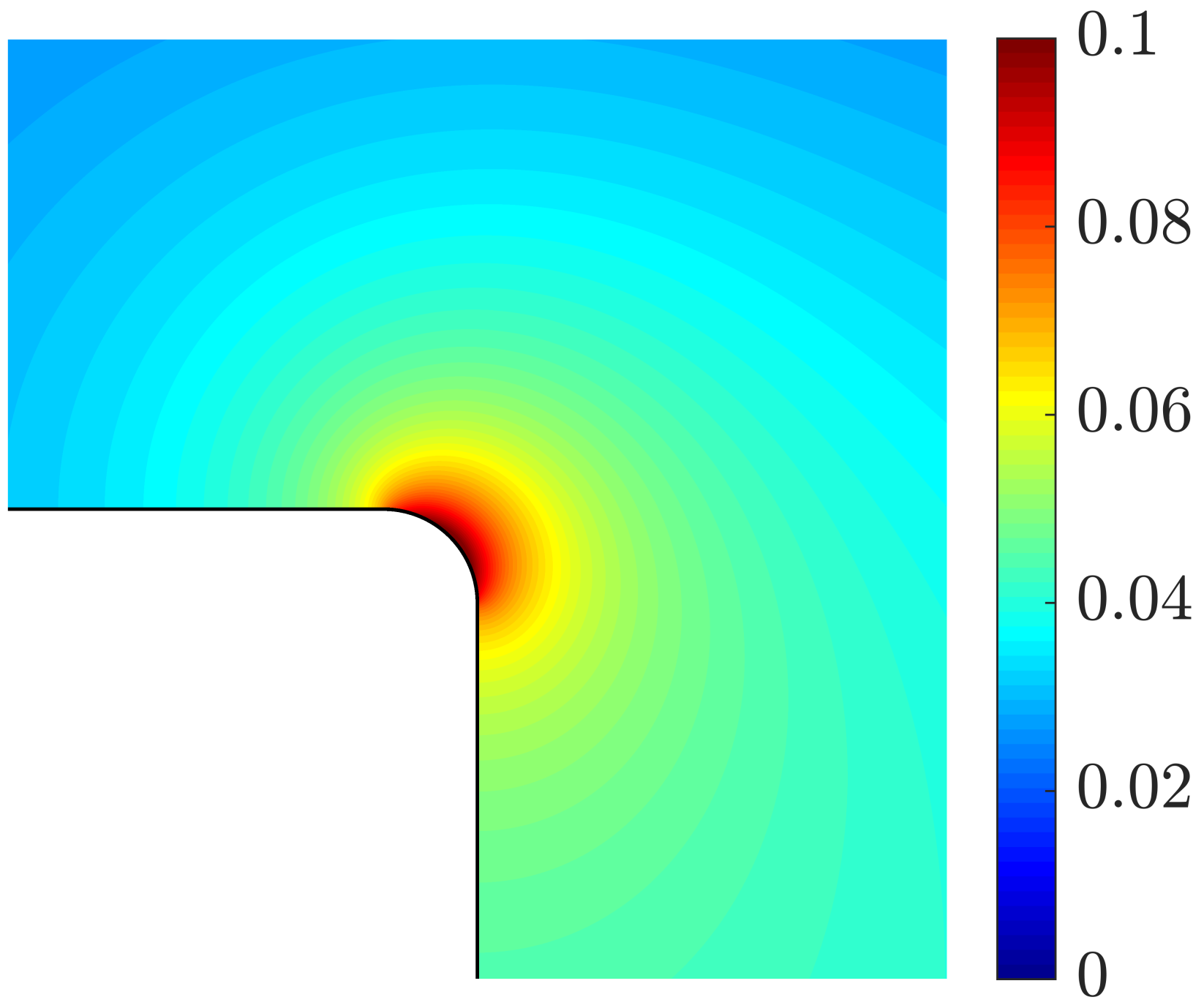}
	\includegraphics[width=0.3\textwidth]{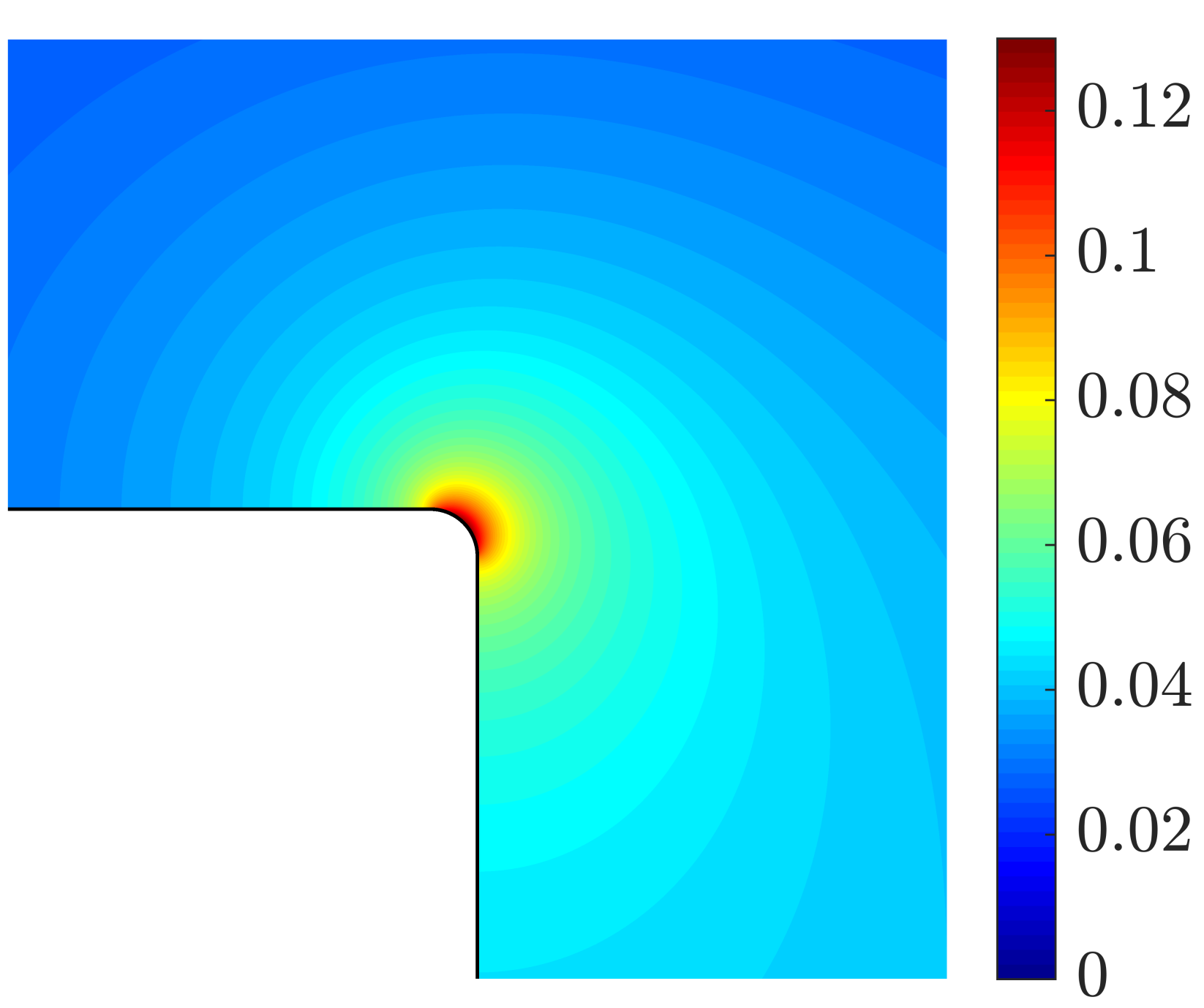}
	\caption{Top: intensity of the electric field for three different radius of the fillet. From left to right: $r{=}5$, $r{=}2$ and $r{=}1$.
				  Bottom: detail of the top-right corner of the inclusion and the computed electric field.
	}
	\label{fig:solutionRadius}
\end{figure}

The proposed degree adaptive process is tested by computing the electric field in a rectangular domain with a square inclusion, $\Omega = [-75,75]{\times}[-100,100] {\setminus} [-50,50]^2$. 
A unit potential is imposed on the outer boundary and a zero potential on the inclusion. 
As it is common in practical engineering applications, the corners of the inclusion are rounded to eliminate the singularity induced by the re-entrant corners~\cite{Krahenbuhl-11}.
Specifically, a fillet defined using a small radius $r$ is introduced to increase the regularity of the boundary. Figure~\ref{fig:solutionRadius} shows the intensity of the electric field for three different geometries, with a fillet of radius $r{=}5$, $r{=}2$ and $r{=}1$ respectively.
The results clearly illustrate the change in the maximum intensity of the electric field (Fig.~\ref{fig:solutionRadius} - top), as well as the localized variations at the corners in terms of the radius of the fillet (Fig.~\ref{fig:solutionRadius} - bottom).

\begin{figure}
	\centering
	\includegraphics[width=0.3\textwidth]{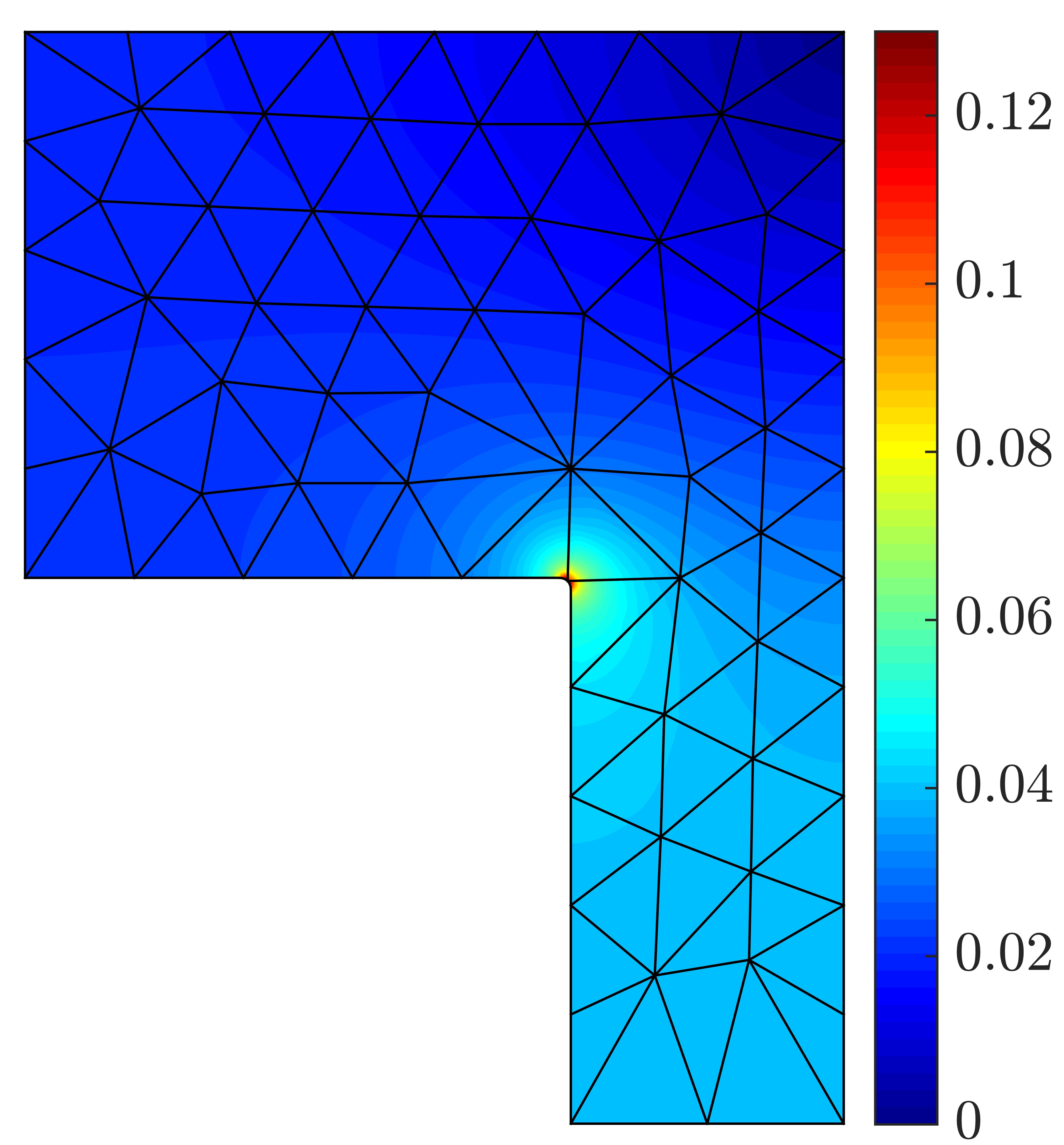}
	\includegraphics[width=0.3\textwidth]{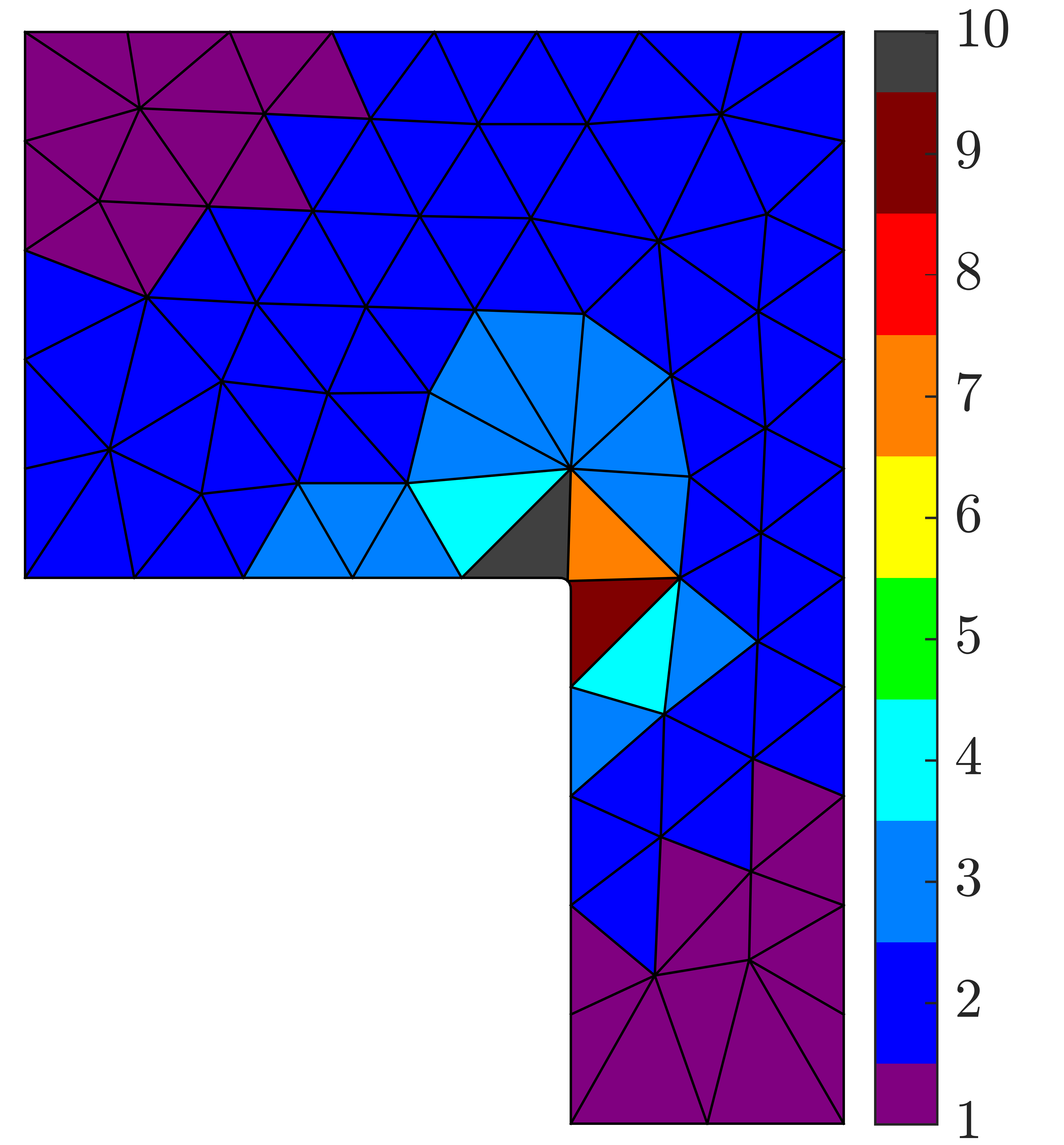}
	\caption{Left: intensity of the electric field computed with the proposed HDG-NEFEM approach. Right: distribution of the approximation degree after eight iterations of the degree adaptive process.}
	\label{fig:solutionAndDegree}
\end{figure}

\begin{figure}
	\centering
	\includegraphics[width=0.48\textwidth]{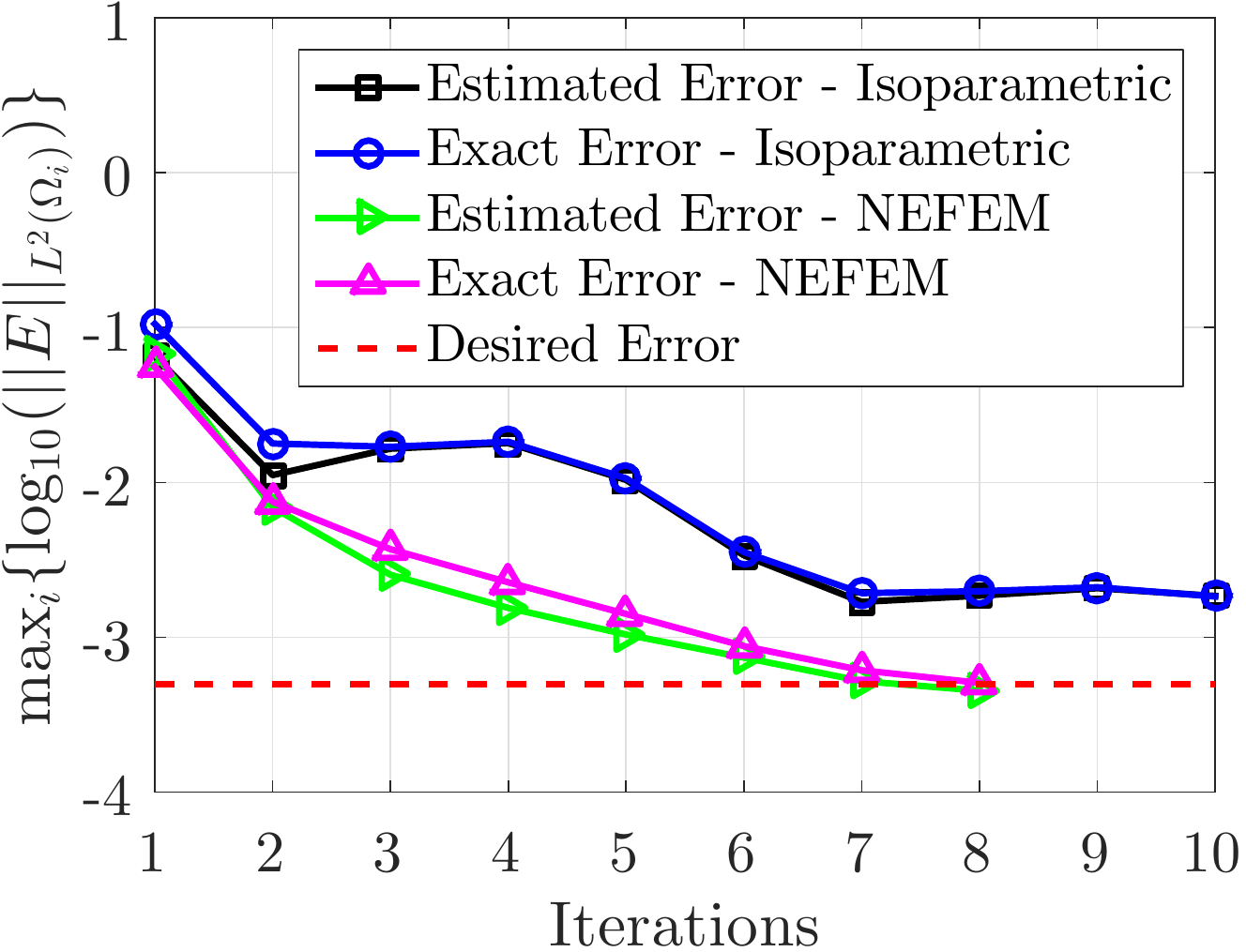}
	\includegraphics[width=0.48\textwidth]{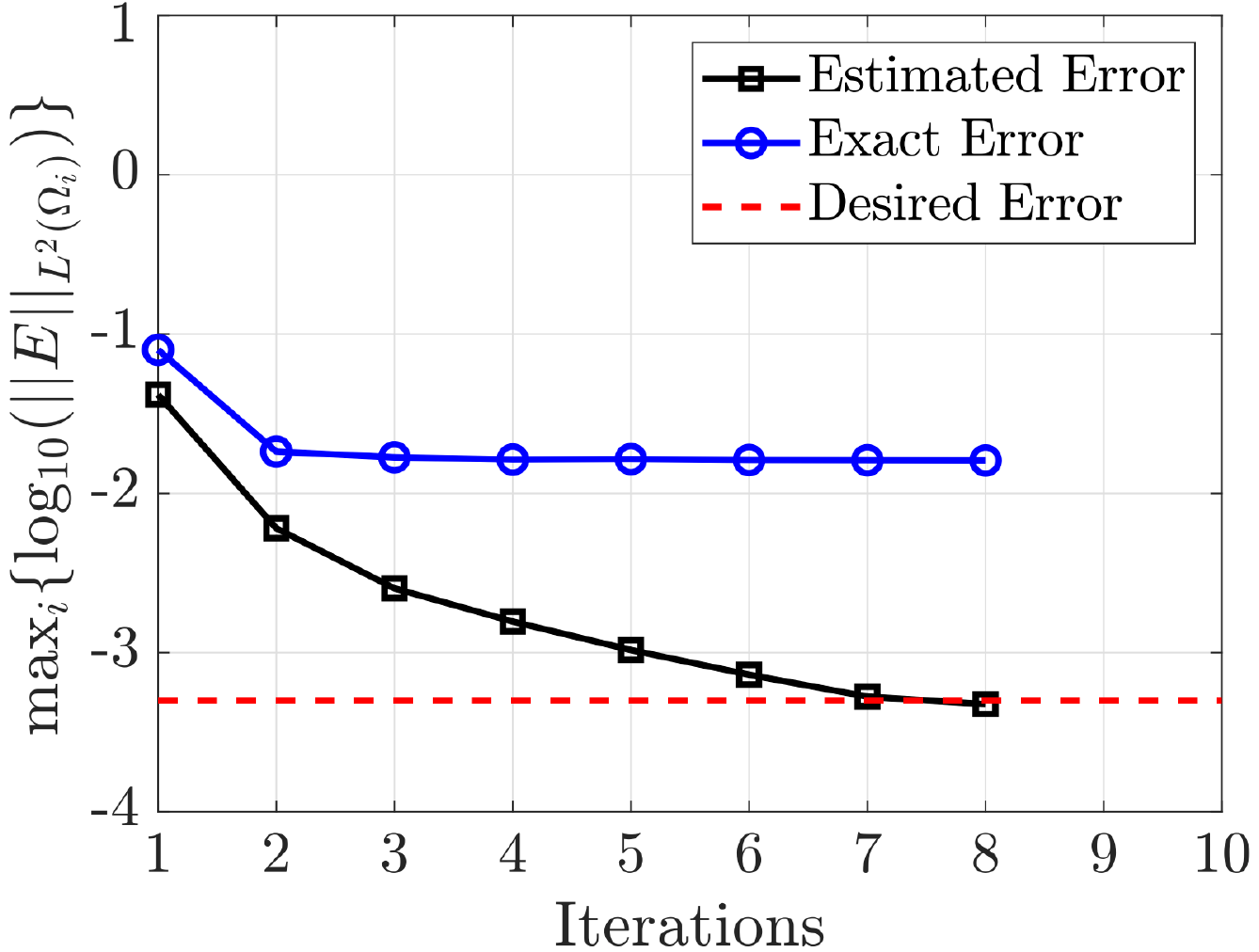}
	
	\caption{Evolution of the exact and estimated errors as a function of the number of iterations in the degree adaptive procedure. Left: isoparametric and NEFEM elements. Right: high-order elements with a fixed approximation of the geometry using cubic polynomials.}
	\label{fig:adaptivityComparison}
\end{figure}

For the application of interest, a fillet of radius $r{=}1$ is considered. In this case, a fine mesh is thus required by isoparametric elements to capture the localized high curvature of the boundary around the corners and the degree adaptive process has to be coupled with mesh adaptation, leading to an $hp$-refinement strategy.
With the proposed HDG-NEFEM approach, a coarse mesh of uniform element size is employed while preserving the exact representation of rounded corners. The degree adaptive process thus determines the required nonuniform degree of approximation to compute the solution with a desired tolerance $\varepsilon {=} 0.5 {\times}10^{-3}$, provided by the user \emph{a priori} and represented by the target dashed line in Figure~\ref{fig:adaptivityComparison}. 
The resulting intensity of the electric field computed on a quarter of the domain and the distribution of the polynomial degree of approximation are depicted in Figure~\ref{fig:solutionAndDegree}.

The evolution of the estimated and exact errors for the proposed HDG-NEFEM approach is shown in Figure~\ref{fig:adaptivityComparison} (left) and compared against the estimated and exact errors for an isoparametric approach. 
It is important to note that the usual isoparametric strategy presents two major drawbacks. 
First, each iteration involving a modification of the polynomial degree in the elements with the rounded corner requires communication with the CAD model and regeneration of the distributions of nodes for the curved elements. 
Second, the change in geometry induced by the change in the degree of the functional approximation is not able to decrease the error towards the imposed tolerance.
The nonsmooth representation of the geometry, i.e. only $\mathcal{C}^0$ between the elements, entails that the numerical approximation of $u$ presents a nonphysical singularity and the degree adaptive process does not provide an optimal solution. 
In contrast, with the proposed HDG-NEFEM approach the error decreases monotonically until the desired tolerance is achieved.
%

To further illustrate the benefits of the proposed HDG-NEFEM approach, a degree adaptive process is performed next with standard high-order elements avoiding the costly communication with the CAD model. To this end, the geometry is represented with cubic polynomials and during the degree adaptive process only the degree of the functional approximation of the solution is changed. 

Figure~\ref{fig:adaptivityComparison} (right) shows the evolution of the estimated and exact errors.
%
%
The results clearly show that despite the adaptive process stops in eight iterations because the estimated error has reached the desired tolerance, the exact error is far from being close to the desired error. This indicates that the adaptive process is actually converging to the solution of a different problem where the geometry is represented with polynomials and remains unchanged during the adaptive iterations. It is worth noticing that using a fixed polynomial approximation of the geometry in a degree adaptive context has been extensively utilized in the literature, see~\cite{giorgiani2013hybridizable,giorgiani2014hybridizable}, but this simple example demonstrates the limitations of such approach.

\section{HDG-Voigt formulation in continuum mechanics and local error indicators}
\label{sc:HDG-Voigt}

In continuum mechanics, the strong enforcement of the symmetry of the stress tensor is associated with the pointwise fulfillment of the conservation of angular momentum.
It is well-known that the classical HDG mixed formulation suffers from suboptimal convergence when low-order discretizations of symmetric second-order tensors are involved~\cite{soon2009hybridizable,Nguyen-CNP:10}.
To remedy this issue, several techniques have been proposed in the context of hybrid discretization techniques~\cite{Ern-DPE-15,Shi-QS-16,Shi-QSS-18,Cockburn-CFQ-17,Cockburn-CF-17}.
The HDG-Voigt formulation introduced in~\cite{RS-SGKH:18,MG-GKSH:18} exploits Voigt notation for second-order tensors, see~\cite{FishBelytschko2007}, to strongly enforce symmetry by storing solely $\msd {:=} \nsd(\nsd {-} 1)/2$ non-redundant off-diagonal components of the stress tensor in a vector form $\stressV$, namely, 
\begin{equation} \label{eq:stressVoigt}
\stressV :=\begin{cases}
\bigl[\sigma_{11} ,\; \sigma_{22} ,\; \sigma_{12} \bigr]^T
& \nsd=2 , \\
\bigl[\sigma_{11} ,\; \sigma_{22} ,\; \sigma_{33} ,\; \sigma_{12} ,\; \sigma_{13} ,\; \sigma_{23} \bigr]^T
& \nsd=3 ,
\end{cases}
\end{equation}
where $\nsd$ is the number of spatial dimensions of the problem.

Consider a domain $\Omega \subset \RR^{\nsd}$ such that $\partial\Omega {=} \Gamma_D \cup \Gamma_N$ and $\Gamma_D \cap \Gamma_N {=} \emptyset$ and the following system of equations describing the behavior of a continuum medium
\begin{equation} \label{eq:Elasticity}
\left\{\begin{aligned}
-\gradS^T \stressV  &= \bm{s}          &&\text{in $\Omega$,}\\
\stressV  &=  \bD(E,\nu) \gradS \bu         &&\text{in $\Omega$,}\\
\bu &= \bu_D     &&\text{on $\Gamma_D$,}\\
\bN^T \stressV &= \bm{g}         &&\text{on $\Gamma_N$,}
\end{aligned} \right.
\end{equation}
where $\bu$ is the unknown displacement field, $\bm{s}$ is the external body force and $\bu_D,\bm{g}$ are the imposed displacement and traction on the boundary, respectively.
The $\msd \times \nsd$ matrices $\gradS$ and $\bN$ account for the linearized symmetric gradient operator and the normal direction to the boundary and have the following form
\begin{subequations}
\begin{equation} \label{eq:symmGrad}
\gradS :=\begin{cases}
\begin{bmatrix}
\partial/\partial x_1 & 0 & \partial/\partial x_2 \\
0 & \partial/\partial x_2 & \partial/\partial x_1
\end{bmatrix}^T
& \nsd=2 , \\
\begin{bmatrix}
\partial/\partial x_1 & 0 & 0 & \partial/\partial x_2 & \partial/\partial x_3 & 0 \\
0 & \partial/\partial x_2 & 0 & \partial/\partial x_1 & 0 & \partial/\partial x_3 \\
0 & 0 & \partial/\partial x_3 & 0 & \partial/\partial x_1 & \partial/\partial x_2
\end{bmatrix}^T
& \nsd=3 . \\
\end{cases}
\end{equation}
\begin{equation} \label{eq:normalVoigt}
\bN :=\begin{cases}
\begin{bmatrix}
n_1 & 0 & n_2 \\
0 & n_2 & n_1
\end{bmatrix}^T
& \nsd=2 , \\
\begin{bmatrix}
n_1 & 0 & 0 & n_2 & n_3 & 0\\
0 & n_2 & 0 & n_1 & 0 & n_3 \\
0 & 0 & n_3 & 0 & n_1 & n_2
\end{bmatrix}^T
& \nsd=3 . \\
\end{cases}
\end{equation}
\end{subequations}
The relationship between the stress tensor and the displacement field is expressed by means of the Hooke's law for linear elastic homogeneous materials, $\stressV = \bD(E,\nu) \gradS \bu$, where the matrix $\bD$ describes the mechanical behavior of the solid as a function of Young's modulus $E$ and Poisson's ratio $\nu$, according to the classical definitions in~\cite{FishBelytschko2007}.

Introducing a symmetric mixed variable $\bL_e$ for the discretization of the strain rate tensor, the linear elasticity problem in Equation~\eqref{eq:Elasticity} is split into a set of $\numel$ local problems that define the primal and mixed variables $(\bu_e,\bL_e)$ as functions of the hybrid variable $\bhu$ representing the trace of the displacement field on the edges/faces of the mesh, namely
\begin{equation} \label{eq:HDG-local}
\left\{\begin{aligned}
\bL_e + \bDHalf \gradS \bu_e &= \bm{0}    &&\text{in $\Omega_e, \, e=1,\ldots ,\numel$,}\\	
\gradS^T \bDHalf \bL_e &= \bm{s}          &&\text{in $\Omega_e, \, e=1,\ldots ,\numel$,}\\
\bu_e &= \bu_D     &&\text{on $\partial\Omega_e \cap \Gamma_D$,}\\
\bu_e &= \bhu     &&\text{on $\partial\Omega_e \setminus \Gamma_D$,}
\end{aligned} \right. 
\end{equation}
and a global problem imposing the Neumann boundary condition and the transmission conditions to enforce inter-element continuity of the solution and the tractions
\begin{equation} \label{eq:HDG-global}
\left\{\begin{aligned}
\bN^T \bDHalf \bL_e &= -\bm{g}         &&\text{on $\Gamma$,}\\
\jump{\bu {\otimes} \bn} &= \bm{0} &&\text{on $\Gamma$,}\\
\jump{\reallywidehat{\bN^T \bDHalf \bL_e}} &= \bm{0}  &&\text{on $\Gamma$,}
\end{aligned} \right. 
\end{equation}
where $\bn$ is the outward normal vector to the faces of the internal skeleton $\Gamma$ and $\reallywidehat{\bN^T \bDHalf \bL_e}$ is the trace of the numerical flux, defined as a function of $\bhu$ and the stabilization parameter $\tau$
\begin{equation} \label{eq:traceL}
\reallywidehat{\bN^T \bDHalf \bL_e} := 
\begin{cases}
\bN^T \bDHalf \bL_e + \tau (\bu_e - \bu_D) & \text{on $\partial\Omega_e\cap\Gamma_D$,} \\
\bN^T \bDHalf \bL_e + \tau (\bu_e - \bhu) & \text{elsewhere.}  
\end{cases}
\end{equation}

Note that owing to the Voigt framework, the mixed variable utilized in the HDG formulation is the pointwise symmetric strain rate tensor, the conservation of angular momentum is fulfilled pointwise in each mesh element and physical tractions are imposed on the Neumann boundary.
Following the HDG rationale, first, the global problem in Equation~\eqref{eq:HDG-global} is solved to obtain $\bhu$ on the internal skeleton $\Gamma$ and on the Neumann boundary $\Gamma_N$. Then, the primal and mixed variables $(\bu_e,\bL_e)$ are computed element-by-element by solving the HDG local problems in Equation~\eqref{eq:HDG-local} independently in each $\Omega_e, \, e{=}1,\ldots,\numel$. Eventually, the following postprocessing procedure is devised: for $e{=}1,\ldots,\numel$ compute a displacement field $\bu_e^\star$ using a polynomial approximation of degree $k{+}1$ such that
\begin{equation} \label{eq:HDGPoissonPostVoigt}
\left\{\begin{aligned}
\gradS^T\bDHalf\gradS\bu_e^\star &= -\gradS^T\bL_e  &&\text{in $\Omega_e$,}\\
\bN^T \bDHalf \gradS \bu_e^\star &= -\bN^T \bL_e  &&\text{on $\partial\Omega_e$,}\\
\end{aligned}\right.
\end{equation}
with the solvability constraint in Equation~\eqref{eq:PostProcessCondMean} to remove the underdetermination due to rigid body translations and 
\begin{equation}\label{eq:PostProcessCondRot}
\int_{\Omega_e}{\!\!\! \grad {\times} \bu_e^\star \, d\Omega} = \int_{\partial\Omega_e}{\!\!\! \bu_e {\cdot} \bt \, d\Gamma} ,
\end{equation}
to account for rigid body rotations, where $\bt$ is the tangential direction to the boundary $\partial\Omega_e$.

It is worth recalling that HHO and HDG, that is, both primal and mixed formulations of hybrid discretization methods display a robust behavior for nearly incompressible materials and do not experience locking phenomena~\cite{Ern-DPE-15,soon2009hybridizable,Fu-FCS-15}. The discussed HDG-Voigt strategy inherits such property. 
Nonetheless, classical HDG methods using approximations with equal-order polynomials of degree $k$ for all the variables experience suboptimal behavior for $k{<}3$.
On the contrary, the proposed HDG-Voigt formulation provides a discretization with optimal convergence of order $k{+}1$ for $\bu,\bL$ and $\bhu$, even in case of low-order polynomial approximations.
Thus, in this context, the advantages of using the HDG-Voigt formulation are twofold. On the one hand, an approximation of the strain rate tensor is directly obtained from the mixed formulation without the need to postprocess the primal variable of the problem. On the other hand, the resulting method provides optimally convergent stress and superconvergent displacement field using a nodal-based approximation for all the variables~\cite{RS-SGKH:18} and without resorting to different interpolation degrees~\cite{Shi-QSS-18} or to the enrichment of the local discrete spaces discussed in~\cite{Cockburn-CF-17}.

\begin{figure}
	\centering
	\includegraphics[width=0.46\textwidth]{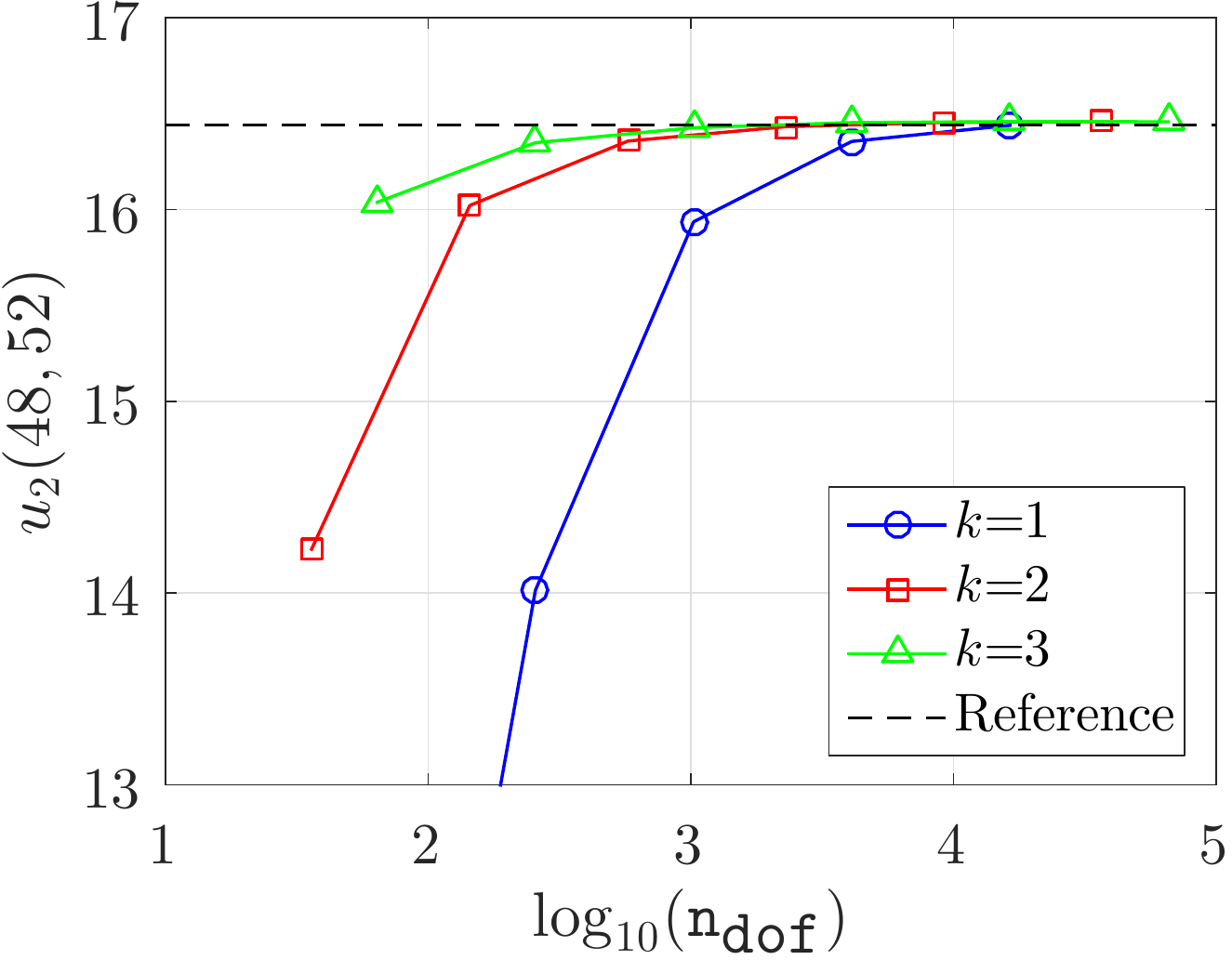}
	\includegraphics[width=0.46\textwidth]{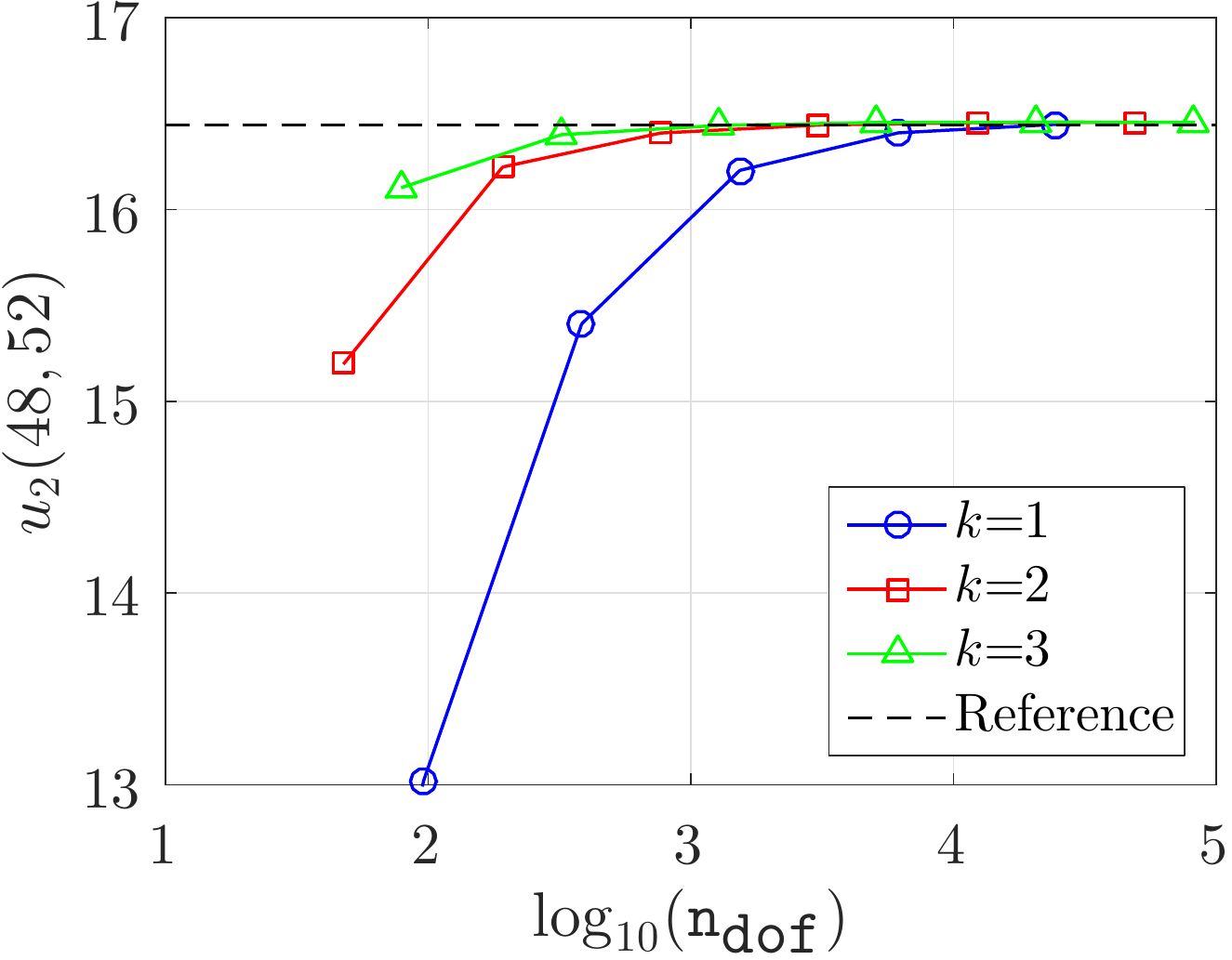}
	\caption{Convergence of the displacement of the mid-point of the right end of Cook's membrane as a function of the number of degrees of freedom of the HDG discretization, using polynomial approximations of degree $k{=}1,2,3$. Left: quadrilateral elements. Right: triangular elements.}
	\label{fig:cooksMembraneDispTip}
\end{figure}

The HDG-Voigt formulation is tested on a well-known benchmark test for bending-dominated elastic problems, the Cook's membrane~\cite{cook2001concepts}.
The domain consists of a tapered plate clamped on the left end and subject to a vertical shear load $\bm{g} {=} (0, 1/16)$ on the opposite end, whereas zero tractions are imposed on the top and bottom parts of the boundary.
Following the problem setup in~\cite{auricchio2005analysis}, a nearly incompressible material with Young modulus $E {=} 1.12499998125$ and Poisson ratio $\nu {=} 0.499999975$ is considered. Figure~\ref{fig:cooksMembraneDispTip} shows the displacement of the mid-point of the right end of the membrane for linear, quadratic and cubic elements on both quadrilateral and triangular meshes. 
The results display the convergence to the reference value, taken from~\cite{auricchio2005analysis}, even for low-order triangular elements, showing the robustness of the HDG-Voigt formulation in the incompressible limit. 

Exploiting the optimal convergence of order $k {+} 1$ of the discretized strain rate tensor $\bL_e$ and the postprocessing procedure discussed in~\cite{RS-SGKH:18,MG-GKSH:18} to resolve the underdetermination due to rigid body motions, a superconvergent approximation $\bu_e^\star$ of the displacement field is constructed.
Thus, the error indicator in Equation~\eqref{eq:errorMeasureU} is computed starting from the approximated primal and postprocessed displacement fields.
Alternatively, a local error indicator based on the strain rate tensor 
\begin{equation} \label{eq:errorMeasureL}
E_e^L = \left[ \frac{1}{|\Omega_e|} \int_{\Omega_e} \left( \gradS \bu_e^\star - \gradS \bu_e \right) {\cdot} \left( \gradS \bu_e^\star - \gradS \bu_e \right) d\Omega \right]^{1/2}
\end{equation}
can be used when a certain level of accuracy is required on the stress tensor rather than on the displacement field~\cite{RS-19}.

\begin{figure}
	\centering
	\includegraphics[width=0.4\textwidth]{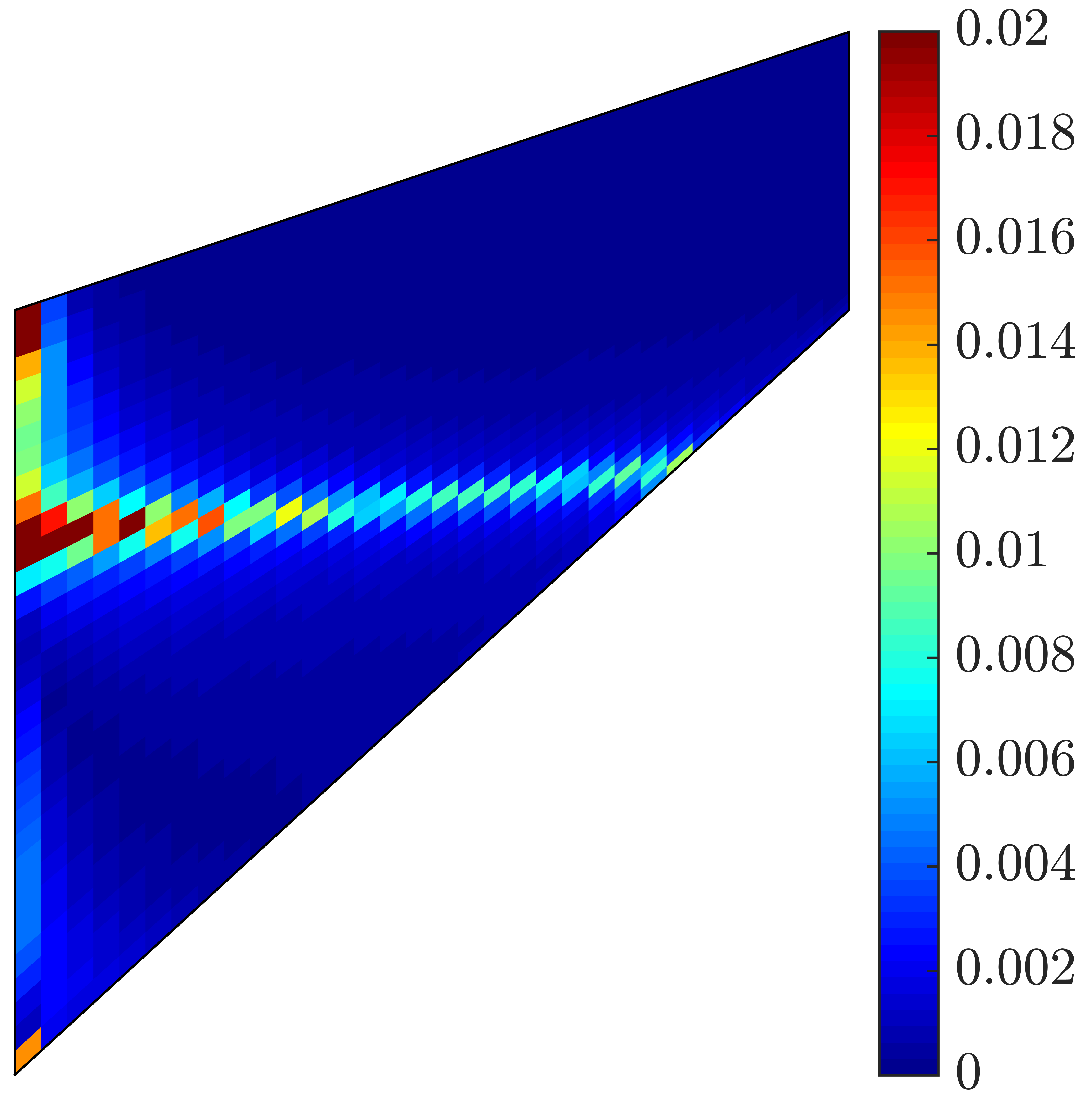}
	\includegraphics[width=0.4\textwidth]{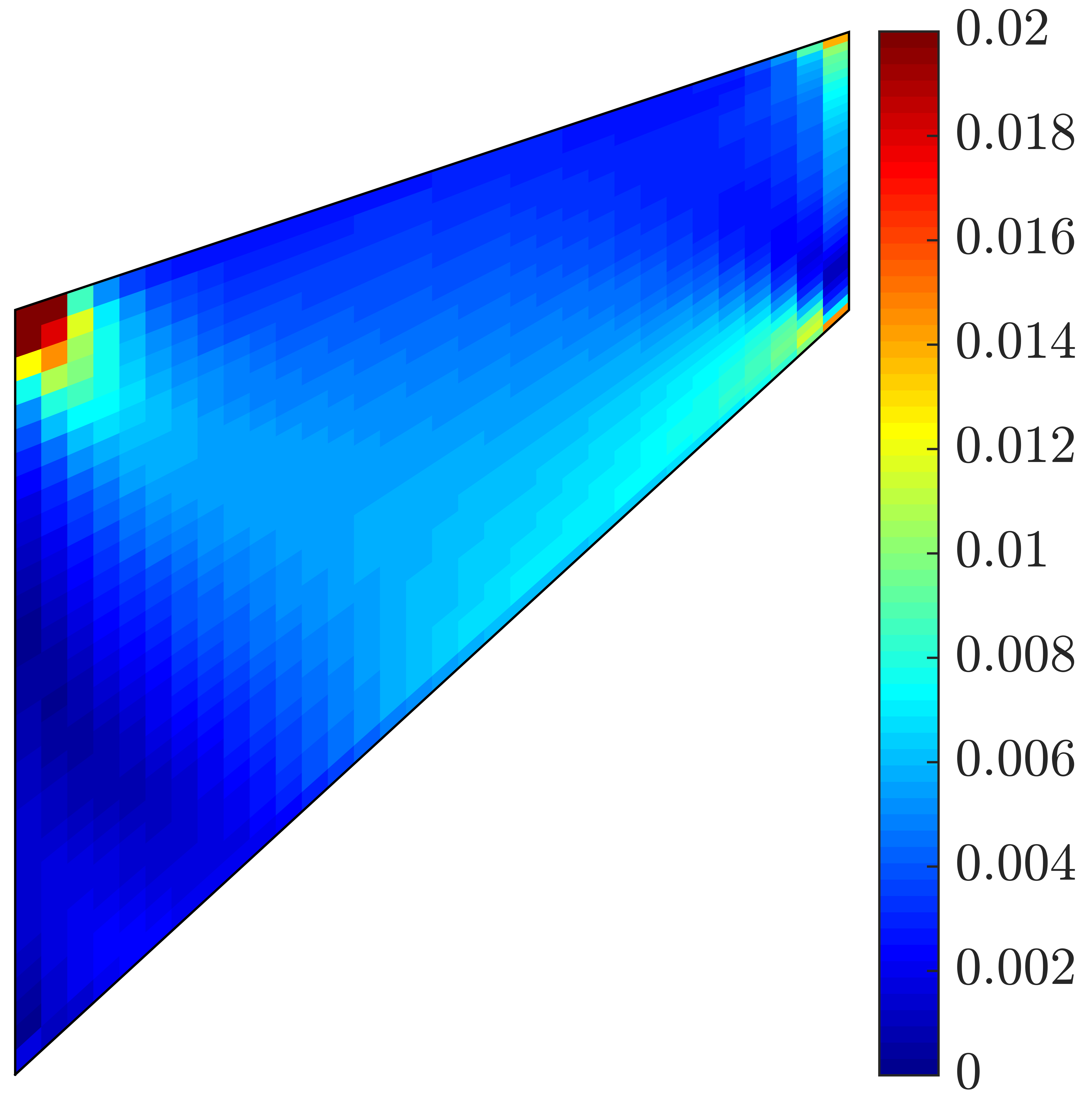}
	\caption{Error indicator based on the displacement field (left) and the strain rate tensor (right).}
	\label{fig:cookQUAH5P1_Indicator}
\end{figure}

Figure~\ref{fig:cookQUAH5P1_Indicator} shows a comparison of the error indicators~\eqref{eq:errorMeasureU} and~\eqref{eq:errorMeasureL} for the displacement field and the strain rate tensor, respectively. 
The different information captured by each error indicator is clearly observed. 
In particular, it is straightforward to observe that the error indicator based on the strain rate tensor is able to provide information about regions where a concentration of stress is present. 
This information is of great interest in engineering applications, e.g. for the optimal design of elastic structures~\cite{Allaire-AD-14}.

\section{FCFV: lowest-order HDG method for large-scale problems}
\label{sc:FCFV}

One of the major challenges that current techniques in computational mechanics face when confronted to industrial applications is proving their ability to efficiently solve large-scale problems in a reliable and robust way.
Despite the numerous advantages in terms of accuracy, efficient treatment of convection-dominated phenomena in flow problems and flexibility for parallelization, the adoption of high-order methods by the industry is still limited, partially due to the difficulty to generate high-order curvilinear meshes of complex configurations~\cite{RS-XSHM-13}.
Starting from the framework discussed above, a novel efficient and robust finite volume (FV) rationale has been proposed in~\cite{RS-SGH:18,RS-SGH:19}.

In order to describe this approach, an incompressible Stokes flow is considered
\begin{equation} \label{eq:Stokes}
\left\{\begin{aligned}
- \grad {\cdot} (\nu \grad \bu -  p \Insd) &= \bm{s}       &&\text{in $\Omega$,}\\
\grad {\cdot} \bu &= 0  &&\text{in $\Omega$,}\\
\bu &= \bu_D  &&\text{on $\Gamma_D$,}\\
\bigl(\nu \grad \bu  -  p \Insd \bigr) \bn &= \bm{g}  &&\text{on $\Gamma_N$,}\\                                          
\end{aligned}\right.
\end{equation}
where the pair $(\bu,p)$ represents the unknown velocity and pressure fields, $\nu {>} 0$ is the viscosity of the fluid, $\Insd$ is the $\nsd {\times} \nsd$ identity matrix and $\bm{s}, \bu_D, \bm{g}$ respectively are the source term, the imposed velocity and pseudo-tractions, see~\cite{donea2003finite}, on the boundary.

Following the HDG rationale introduced in Section~\ref{sc:HDG}, the FCFV local and global problems for the Stokes equations are introduced. More precisely, in each cell $\Omega_e, \ e=1,\ldots,\numel$, it holds
\begin{equation} \label{eq:FCFV-local}
\left\{\begin{aligned}
\bL_e + \sqrt{\nu} \grad \bu_e &= \bm{0}    &&\text{in $\Omega_e, \,e=1,\ldots ,\numel$,}\\	
\grad {\cdot} \left(\sqrt{\nu} \bL_e + p_e \Insd \right) &= \bm{s}          &&\text{in $\Omega_e, \,e=1,\ldots ,\numel$,}\\
\grad {\cdot} \bu_e &= 0         &&\text{in $\Omega_e, \, e=1,\ldots ,\numel$,}\\
\bu_e &= \bu_D     &&\text{on $\partial\Omega_e \cap \Gamma_D$,}\\
\bu_e &= \bhu     &&\text{on $\partial\Omega_e \setminus \Gamma_D$,}
\end{aligned} \right. 
\end{equation}
with the following additional constraint to remove the underdetermination of pressure due to the Dirichlet boundary conditions imposed in Equation~\eqref{eq:FCFV-local}
\begin{equation} \label{eq:constraint-local}
\int_{\partial\Omega_e}{\!\!\! p_e \, d\Gamma} = |\partial\Omega_e| \rho_e .
\end{equation}
The FCFV global problem features the Neumann boundary conditions and the transmission conditions enforcing inter-element continuity of the solution and the fluxes, as previously detailed for the HDG method
\begin{equation} \label{eq:FCFV-global}
\left\{\begin{aligned}
\left( \sqrt{\nu}\bL_e + p_e\Insd \right) \! \bn &= -\bm{g}         &&\text{on $\Gamma$,}\\
\jump{\bu {\otimes} \bn} &= \bm{0} &&\text{on $\Gamma$,}\\
\jump{\reallywidehat{\left( \sqrt{\nu}\bL_e + p_e\Insd \right) \! \bn}} &= \bm{0}  &&\text{on $\Gamma$,}
\end{aligned} \right. 
\end{equation}
where the numerical normal flux on the boundary is defined as
\begin{equation} \label{eq:traceStokes}
\reallywidehat{\left( \sqrt{\nu}\bL_e + p_e\Insd \right) \! \bn} := 
\begin{cases}
\left( \sqrt{\nu} \bL_e + p_e \Insd \right) \! \bn + \tau (\bu_e - \bu_D) & \text{on $\partial\Omega_e\cap\Gamma_D$,} \\
\left( \sqrt{\nu} \bL_e + p_e \Insd \right) \! \bn + \tau (\bu_e - \bhu) & \text{elsewhere.}  
\end{cases}
\end{equation}
Moreover, the incompressibility constraint is expressed in weak form as
\begin{equation} \label{eq:constraint-global}
\int_{\partial\Omega_e \setminus \Gamma_D}{\!\!\! \bhu {\cdot} \bn \, d\Gamma} + \int_{\partial\Omega_e \cap \Gamma_D}{\!\!\! \bu_D {\cdot} \bn \, d\Gamma} = 0 \qquad \text{for $e=1,\ldots ,\numel$}.
\end{equation}

FCFV may be interpreted as the lowest-order HDG mixed method which employs a constant degree of approximation in each cell for the velocity $\bu_e$, the pressure $p_e$ and the mixed variable $\bL_e$, representing the gradient of velocity, a constant degree of approximation on each edge/face for the velocity $\bhu$ and a constant value $\rho_e$ for the mean pressure in each cell. 
Moreover, the FCFV global and local problems are discretized using a quadrature with one integration point located in the centroid of the cell or face and in the midpoint of the edge.

FCFV solves the problem in two phases~\cite{RS-SGH:18,RS-SGH:19}.
First, by applying the divergence theorem to Equation~\eqref{eq:FCFV-local} and exploiting the definition of the numerical normal flux on the boundary in Equation~\eqref{eq:traceStokes}, a set of $\numel$ local integral problems is obtained
\begin{equation} \label{eq:integral-local}
\begin{aligned}
	- \int_{\Omega_e}{\!\!\! \bL_e \, d\Omega} &=   \int_{\partial\Omega_e\cap\Gamma_D}{\!\!\! \sqrt{\nu}\bu_D {\otimes} \bn \, d\Gamma} + \int_{\partial\Omega_e\setminus\Gamma_D}{\!\!\! \sqrt{\nu} \bhu {\otimes} \bn \, d\Gamma} ,
	\\
	\int_{\partial\Omega_e}{\!\!\! \tau \, \bu_e \, d\Gamma} &=  \int_{\Omega_e}{\!\!\! \bm{s} \, d\Omega} + \int_{\partial\Omega_e\cap\Gamma_D}{\!\!\! \tau \, \bu_D \, d\Gamma} + \int_{\partial\Omega_e\setminus\Gamma_D}{\!\!\! \tau \, \bhu \, d\Gamma} ,
	\\
	 \int_{\partial\Omega_e}{\!\!\! p_e \, d\Gamma} &= |\partial\Omega_e| \rho_e .
\end{aligned}
\end{equation}
Note that the divergence theorem applied to the incompressibility constraint in Equation~\eqref{eq:FCFV-local} leads to Equation~\eqref{eq:constraint-global}, which is thus omitted from the local problem since only the global unknown $\bhu$ is involved. Moreover, the last equation of the previous system directly stems from Equation~\eqref{eq:constraint-local}.

It is worth noticing that the equations of the FCFV local problem decouple and a closed-form expression of all the variables as functions of the velocity $\bhu$ on the boundary $\partial\Omega_e {\setminus} \Gamma_D$ and the mean value $\rho_e$ of the pressure inside the element $\Omega_e$ is obtained.
The previously determined elemental expressions of $(\bu_e,p_e,\bL_e) $ are employed to solve the FCFV global problem~\eqref{eq:FCFV-global} with the incompressibility constraint in Equation~\eqref{eq:constraint-global}, namely
\begin{equation} \label{eq:integral-global}
\begin{aligned}
& \sum_{e=1}^{\numel}\left\{
	\int_{\partial\Omega_e\setminus\Gamma_D}{\!\!\! \sqrt{\nu} \bL_e \bn \, d\Gamma}	
	+ \int_{\partial\Omega_e\setminus\Gamma_D}{\!\!\! p_e \bn \, d\Gamma}		
	+ \int_{\partial\Omega_e\setminus\Gamma_D}{\!\!\! \tau \, \bu_e \, d\Gamma} \right.
	\\
& \left. \hspace{120pt}	
- \int_{\partial\Omega_e\setminus\Gamma_D}{\!\!\! \tau \, \bhu \, d\Gamma} \right\}
	= -\sum_{e=1}^{\numel} \int_{\partial\Omega_e\cap\Gamma_N}{\!\!\! \bm{g} \, d\Gamma} ,
	\\
&
	\int_{\partial\Omega_e \setminus \Gamma_D}{\!\!\! \bhu {\cdot} \bn \, d\Gamma} = - \int_{\partial\Omega_e \cap \Gamma_D}{\!\!\! \bu_D {\cdot} \bn \, d\Gamma} = 0 \qquad \text{for $e=1,\ldots ,\numel$}.
\end{aligned}
\end{equation}

The resulting linear system obtained from the FCFV discretization of Equation~\eqref{eq:integral-global} is symmetric and features a saddle-point structure with $\nface\nsd {+} \numel$ unknowns, being $\nface$ the number of internal and Neumann edges/faces.
The FCFV global problem has a sparse block structure allowing a computationally efficient implementation, see~\cite{RS-SGH:18}. Moreover, FCFV local computations to determine velocity, pressure and gradient of velocity in the centroid of each cell solely involve elementary operations cell-by-cell for which modern parallel architectures can be exploited.

FCFV inherits the approximation properties of the corresponding high-order HDG formulation from which it is derived.
More precisely, optimal first-order convergence is obtained for velocity, pressure and gradient of velocity. In addition, contrary to other mixed finite element methods, with the FCFV it is possible to use the same space of approximation for both velocity and pressure, circumventing the so-called Ladyzhenskaya-Babu{\v s}ka-Brezzi  (LBB) condition.

Compared to other FV methods, the FCFV provides first-order accuracy of the solution and its gradient without the need to perform flux reconstruction  as in the context of cell-centered and vertex-centered finite volumes~\cite{diskin2010comparison,diskin2011comparison}.
Furthermore, the accuracy of the FCFV method is preserved in presence of unstructured meshes, with distorted and stretched cells~\cite{RS-SGH:18,RS-SGH:19}. This is of major importance when solving problems in complex geometries as other FV methods lose accuracy and optimal convergence properties when non-orthogonal and anisotropic cells are introduced in the computational mesh~\cite{diskin2010comparison,diskin2011comparison}.

To highlight the efficiency of the proposed FCFV method, a Stokes flow is simulated in a channel with $39$ rigid particles in the shape of red-blood cells (RBCs). 
A parabolic velocity profile modelling an undisturbed flow is imposed on the inlet and on the outlet of the channel, whereas a no-slip boundary condition is imposed on the remaining walls and on the surface of the particles. 

\begin{figure}
	\centering
	\includegraphics[width=0.68\textwidth]{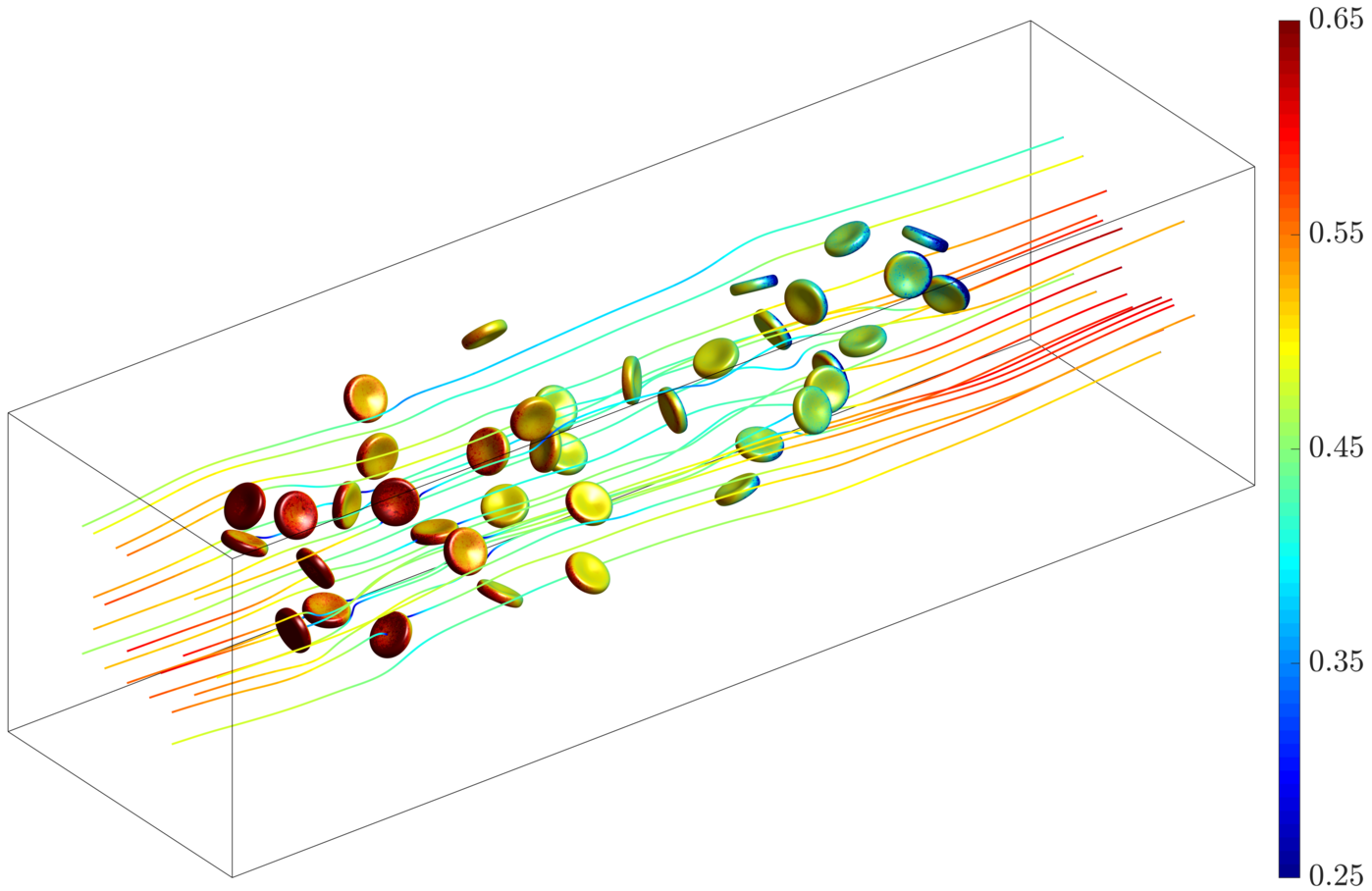}
	
	\caption{Pressure field on $39$ particles modelling RBCs immersed in an incompressible Stokes flow in a channel and streamlines of the velocity field.}
	\label{fig:RBC-pressure}
\end{figure}

The computational domain $\Omega {=} [-10,25] {\times} [-5,5] {\times} [-5,5] {\setminus} \mathcal{B}$, where $\mathcal{B}$ is the union of the $39$ RBCs, is discretized using an unstructured mesh of $8,972,888$ tetrahedral cells, $35,891,552$ nodes and $17,523,981$ internal faces.
The FCFV global system for the mesh configuration under analysis features $61,544,832$ unknowns.
The simulation was performed using a code developed in Matlab$^{\text{\textregistered}}$.
The computation of all the elemental contributions to the global system took 51 minutes whereas 18 minutes were required for the assembly of the matrix.
The solution of the linear system was performed using the Matlab$^{\text{\textregistered}}$ biconjugate gradient method in a single processor and without preconditioner.
%
Eventually, the evaluation of the element-by-element solution in all 61 millions elements took 7 minutes using a single processor.

\begin{figure}
	\centering
	\includegraphics[width=0.68\textwidth]{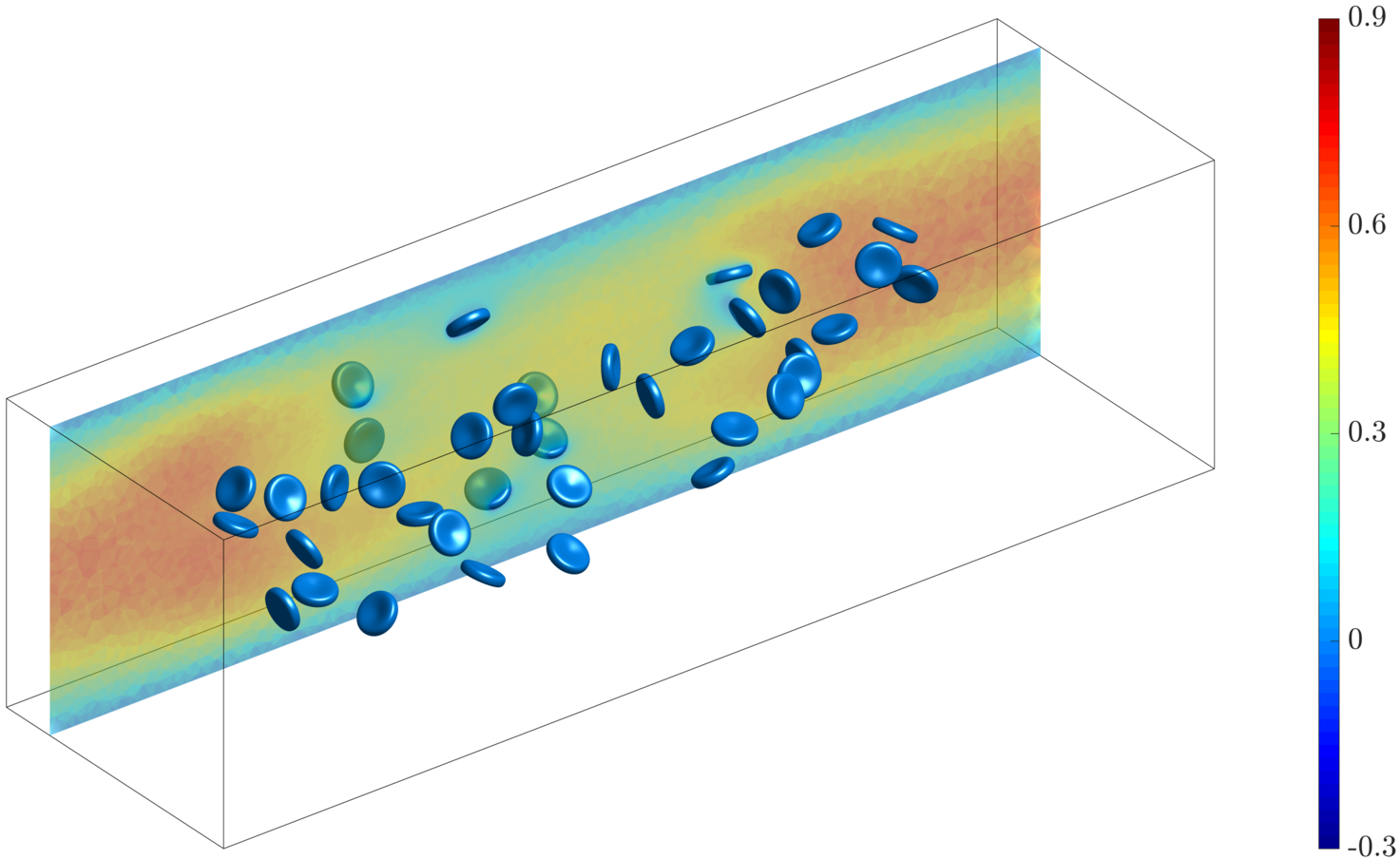}
	
	\includegraphics[width=0.68\textwidth]{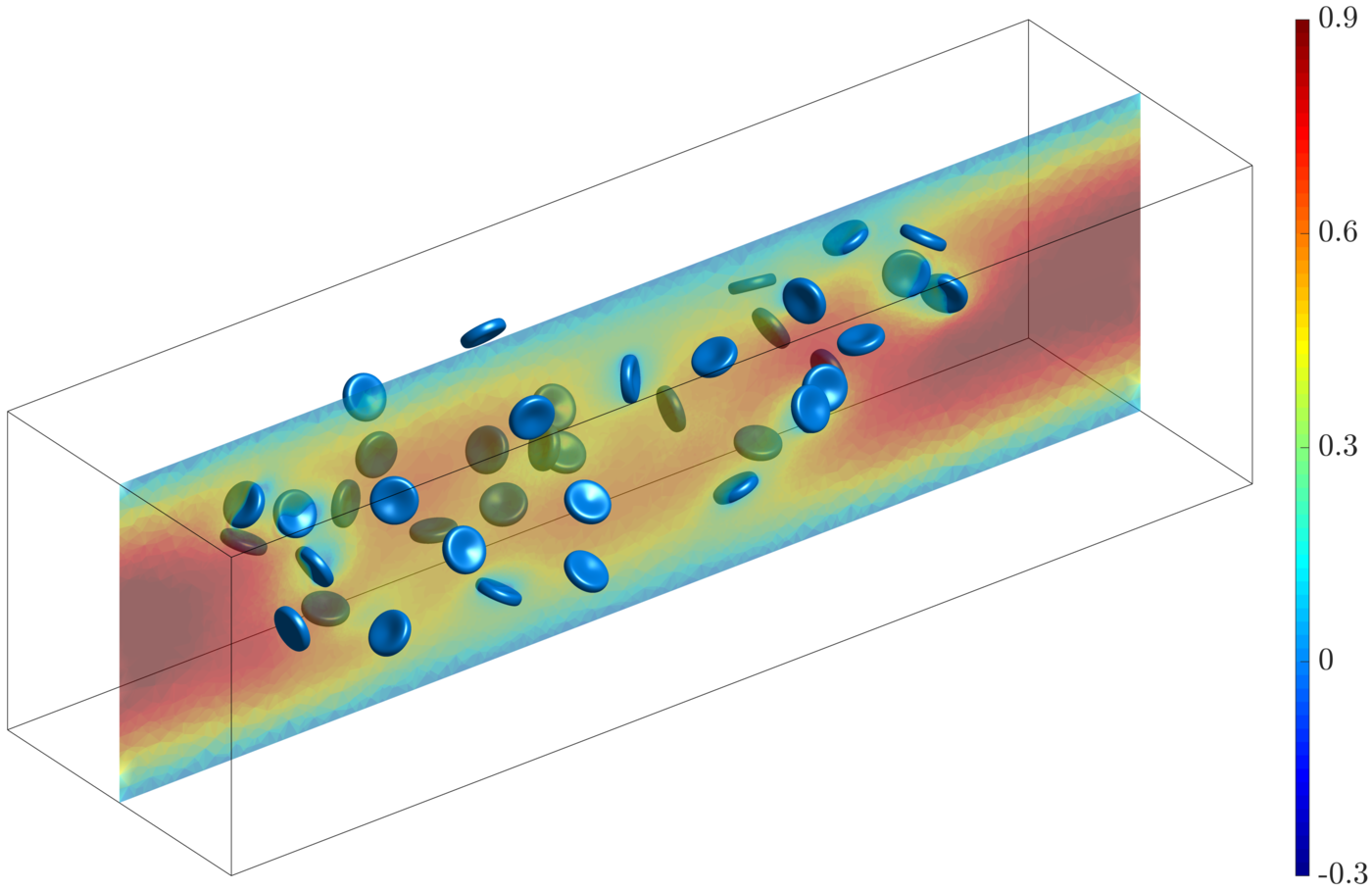}
	
	\includegraphics[width=0.68\textwidth]{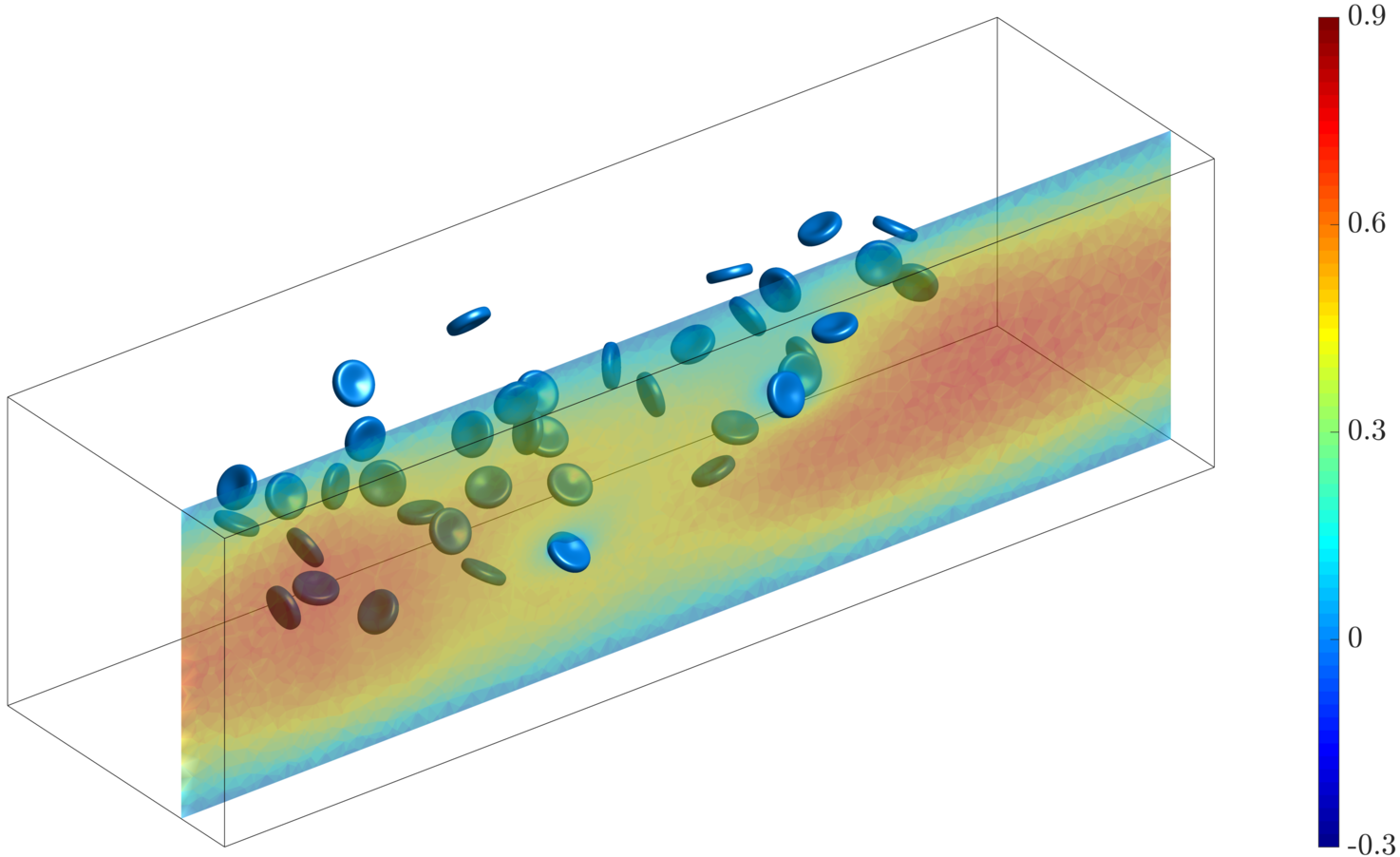}
	
	\caption{Magnitude of the velocity field of an incompressible Stokes flow in a channel with $39$ particles modelling RBCs. Section plane for $y{=}{-}3$ (top), $y{=}0$ (center) and $y{=}3$ (bottom).}
	\label{fig:RBC-velocity}
\end{figure}

The pressure distribution on the surface of the RBCs and the velocity streamlines are presented in Figure~\ref{fig:RBC-pressure}. 
%
Figure~\ref{fig:RBC-velocity} displays the magnitude of the velocity field at three different sections of the computational domain.

\section{Concluding remarks}
\label{sc:Conclusion}

Three recent contributions to HDG are discussed in this paper.
First, the HDG-NEFEM paradigm exploits the description of the boundary of the domain via NURBS to construct an HDG approximation with exact geometry and to devise an efficient and robust degree adaptivity strategy.
Second, the HDG-Voigt formulation is utilized in the context of continuum mechanics to devise an HDG method with pointwise symmetric mixed variable, namely the strain rate tensor. The resulting formulation allows to achieve optimal convergence and superconvergence properties even for low-order polynomial approximations and, consequently, to compute local error indicators based on either the displacement or the stress field.
Third, the FCFV rationale proposes a fast implementation of the lowest-order HDG method. The resulting finite volume paradigm is reconstruction-free, robust to mesh distortion and element stretching and is able to efficiently tackle large-scale problems.
Ongoing investigations focus on the application of the discussed strategies to nonlinear problems of interest in engineering applications and to the simulation of transient phenomena.

\section*{Acknowledgements}
This work is partially supported by the European Union's Horizon 2020 research and innovation programme under the Marie Sk\l odowska-Curie actions (Grant No. 675919) and the Spanish Ministry of Economy and Competitiveness (Grant No. DPI2017-85139-C2-2-R). The first author also gratefully acknowledges the financial support provided by Generalitat de Catalunya (Grant No. 2017-SGR-1278).

\paragraph*{Conflict of interest} \ The authors declare that they have no conflict of interest.

\bibliographystyle{unsrt}
\bibliography{Ref-HDG}

\end{document}